# Homotopy Continuation Enhanced Branch and Bound Algorithms for Strongly Nonconvex Mixed-Integer Nonlinear Programming Problems

Yingjie Ma, and Jie Li[*]

Centre for Process Integration, Department of Chemical Engineering and Analytical Science, School of Engineering, The University of Manchester, Manchester M13 9PL, UK

## Abstract

Large-scale strongly nonlinear and nonconvex mixed-integer nonlinear programming (MINLP) models frequently appear in optimisation-based process synthesis, integration, intensification and process control. However, they are usually difficult to solve by existing algorithms within acceptable time. In this work, we propose two robust homotopy continuation enhanced branch and bound (HCBB) algorithms (denoted as HCBB-FP and HCBB-RB) where the homotopy continuation method is employed to gradually approach the optimum of the NLP subproblem at a node from the solution at its parent node. A variable step length is adapted to effectively balance feasibility and computational efficiency. The computational results from solving four existing process synthesis problems demonstrate that the proposed HCBB algorithms can find the same optimal solution from different initial points, while the existing MINLP algorithms fail or find much worse solutions. In addition, HCBB-RB is superior to HCBB-FP due to much lower computational effort required for the same locally optimal solution.

[*] Corresponding author:  Jie Li (Email address: jie.li-2@manchester.ac.uk). Tel: + 44 (0) 161 306 8622. Fax: +44 (0) 161 236 7439

**Introduction**

Mixed-integer nonlinear programming (MINLP) is one of the most general class of optimisation problems[1], which appear in many applications of engineering, applied mathematics and operations research[2]. Especially, it is widely used in the field of process systems engineering (PSE)[3], such as process synthesis[4,5], planning and scheduling[6-9], process control[10,11], and moledular design[12,13]. At the same time, the PSE applications have spurred many developments in MINLP[14].

There are two classes of MINLP problems, including convex MINLP and nonconvex MINLP[15]. The convex MINLP problems can be solved more easily as the resulting nonlinear programming (NLP) problems after relaxing the integer variables are convex and the convexity of the objective function and the relaxed feasible regions can be exploited. The nonlinear branch and bound (B&B) method[16,17] solves an NLP subproblem at each node of the B&B tree, where part of the binary variables are fixed, while the others are relaxed as continuous variables. Compared with the nonlinear B&B, decomposition-based methods, such as the generalized Bender decomposition (GBD)[18,19] algorithm and the outer approximation (OA) algorithm[20] usually solve less NLP subproblems with the cost of solving additional mixed-integer linear programming (MILP) problems though. As a result, they are possibly more efficient when NLP subproblems are more expensive to solve[20]. Another class of algorithms are generalized disjunctive programming (GDP)-based algorithms, including GDP B&B[21] and logic-based OA[22], which can reduce the scale of the NLP subproblems and circumvent some numerical difficulties in the nonlinear B&B and OA algorithms[23]. Some other algorithms, such as extended cutting plane[24] and LP/NLP B&B are also reported[25].

The optimality and convergence of most of the above algorithms rely on the preliminaries that the relaxed MINLP (i.e., NLP) subproblems are convex and can be solved when they are feasible. However, the nonconvex MINLP problems, where the NLP subproblems are nonconvex and are even hard to converge to local optima in some cases, are frequently derived to accurately model the real-world problems, particularly in chemical engineering[3,26,27]. Deterministic global optimisation algorithms[26] have been proposed to solve such problems, which are usually based on spatial branch and bound[28-30], where both binary and continuous variables are branched. Although the deterministic algorithms can theoretically provide rigorously global optimal solutions, they may generate a huge search tree for large-scale strongly nonconvex problems, leading to unacceptable computational time to reach $\varepsilon$-global optimality[15]. Due to that, stochastic algorithms, such as simulated annealing[31], genetic algorithms[32], and particle swarm optimisation[33], are frequently used to generate a good solution within affordable time for large-scale strongly nonconvex MINLP problems. However, the stochastic algorithms lead to long computational times and no guarantee of optimality[34]. Another class of methods are to relax an MINLP problem as an NLP problem and add some special constraints to force the integrality of the relaxed integer variables gradually[35]. Such method has been applied to optimistion of distillation columns[36] and synthesis of reaction-separation-recycle processes[37]. However, due to the forcing constraints, the interior of the relaxed integer variables are ill-defined, causing local optimum or even infeasibility without careful initialization[38].

Due to the above reasons, the algorithms developed for the convex MINLP problems are often used to generate a local optimum for large-scale strongly nonconvex MINLP problelms in many PSE applications, especially when detailed unit operation models are employed[39]. For example, the outer approximation algorithm[15] with some modifications[40-42] has been applied to optimisation of both single distillation column[41,43,44] and distillation sequence[45] modeled by rigorous material balance, phase equilibrium, summation and enthalpy balance (MESH) equations, synthesis of heat exchanger network[4], simultaneous process and heat exchanger network optimisation[46], and simultaneous topology and parameter optimization of a multiple cantilever beam[42] and roller and sliding hydraulic steel gate[47]. The logic-based OA algorithm[22] has also been applied to solve the synthesis of distillation sequence and optimize the whole plant using MESH model for distillation[48-50]. The nonlinear B&B has been applied to integrated desing and control[51,52], simultaneous scheduling and control[53], process synthesis incorporating MESH models for distillation columns[54]. However, the difficulty in solving nonconvex

NLP subproblems may still lead to worse solutions or even no solution at all for the nonconvex MINLP problems[52], especially when some accurate yet complex models are included[55]. Thus, when solving large-scale strongly nonconvex MINLP using local optimisation algorithms, good initial points and/or very careful formulations are required[39,45,48,55], which, however, are usually case-by-case and largely rely on trial and error and the experience of the modellers. To resolve such problem, on the one hand, more robust NLP algorithms/methods are desired. For example, a trust region filter (TRF) method has been proposed[56] and refined[57] for optimisation of chemical processes with high-fidelity models, and for heat exchanger network synthesis with detailed exchanger design[34,55]. However, the convergence of the algorithm may largely rely on the robust convergence and accuate derivative information of the simulation using the original model, which may be unavailiable or unreliable in some cases. In addition, which parts of the high-fidelity model are substituted with what kinds of surrogate models might need expert knowledge to achieve good convergence and high efficiency. Another pathway to improve the convergence of NLP subproblems in an MINLP algorithm is to make full use of solutions generated from the preceding NLP subproblems in a sysmetic way, which has not been emphasized to the best of our knowledge and will be done in this work.

An important application of MINLP in chemical engineering is process synthesis[2], which is to determine the optimal interconnections of processing units as well as the optimal type and design of the units within a process system[58]. Many reports on optimisation-based process synthesis use shortcut unit operation models to get solutions within reasonable time[59-65], but the solutions generated from the shortcut model-based formulation may be infeasible, impractical or suboptimal when detailed deisgn models are considered[55]. Thus, rigorous operation models are desired to be incorporated into process synthesis, which usually leads to large-scale strongly nonconvex MINLP problems[58], posing significant challenges to the existing algorithms as discussed above. Actually, this is the main application field motivating the current work.

In this work, we propose two robust homotopy continuation enhanced B&B (HCBB) algorithms to solve large-scale strongly nonconvex MINLP problems in the aforementioned PSE applications especially in process synthesis. In these algorithms, the homotopy continuation (HC) method is employed to improve the convergence and efficiency of solving the NLP subproblem at each node through constructing a homotopy path using the solution obtained at its parent node. This new homotopy path is constructed through exploiting the similarity of the NLP subproblems generated during B&B, which is completely different from the existing HC method for solving single NLP problem, where the HC path was constructed through modifying the Karush-Kuhn-Tucker (KKT) conditions[66,67]. In our method, during HC, an adaptive variable-step length method is used to balance the convergence and computational efficiency, and three match conditions are proposed to exploit the successful HC step lengths at previous nodes for higher efficiency. While the first HCBB algorithm (denoted as HCBB-FP) solves a feasibility problem with a fixed value for the homotopy variable in each HC step, the other one (called HCBB-RB) solves an optimality problem with a relaxed lower or upper bound of the homotopy variable. Four challenging process synthesis problems from literature are solved to evaluate the capability of the proposed HCBB algorithms. The computational studies demonstrate that the proposed HCBB algorithms solve all examples to the same high-quality local optimum from different initial points, whilst all other existing algorithms including the nonlinear branch and bound algorithm in the GAMS/SBB solver[68], the OA algorithm in the GAMS/DICOPT solver[69], and the SRMINLP method from Ma et al.[37] fail or find much worse locally optimal solutions. In addition, HCBB-RB is able to find almost the same locally optimal solution with much lower computational effort compared to HCBB-FP.

The rest of the paper is organized as follows: The next section defines the MINLP problem to be solved in the current work. Then we show the HCBB algorithms in detail. After that, four process synthesis problems are solved to illustrate the convergence and efficiency of the proposed algorithms with a fair comparison with existing algorithms. Finally, we conclude the work with some useful insights.

**Problem Statement**

The large-scale strongly nonconvex MINLP problem to be solved in the current work is stated as follows:

$$\min_{\mathbf{x},\mathbf{y}} f(\mathbf{x},\mathbf{y}) \qquad\qquad (\text{P})$$

$$s.t.\ \mathbf{h}(\mathbf{x},\mathbf{y}) = 0$$

$$\mathbf{g}(\mathbf{x},\mathbf{y}) \leq 0$$

$$\mathbf{x} \in \mathbb{R}^n, \mathbf{y} \in \mathbb{Z}^m,$$

where $\mathbf{x}$ is a vector of $n$ continuous variables, i.e., $\mathbf{x} = [x_1, x_2, \dots, x_n]^T$, and $\mathbf{y}$ is a vector of $m$ integer variables, i.e., $\mathbf{y} = [y_1, y_2, \dots, y_m]^T$. Here, we restrict $\mathbf{y}$ to binary variables without loss of generalization as any integer variables can be expressed by binary variables[2]. $f$ is a scalar function, whilst $\mathbf{g}$ and $\mathbf{h}$ are vector functions.

The relaxation of the original MINLP problem P is provided below,

$$\min_{\mathbf{x},\mathbf{y}} f(\mathbf{x},\mathbf{y}) \qquad\qquad (\text{RP})$$

$$s.t.\ \mathbf{h}(\mathbf{x},\mathbf{y}) = 0$$

$$\mathbf{g}(\mathbf{x},\mathbf{y}) \leq 0$$

$$\mathbf{x} \in \mathbb{R}^n, \mathbf{y} \in \prod_1^m [0,1]$$

We assume that $f(\cdot), \mathbf{g}(\cdot), \mathbf{h}(\cdot)$ are twice differentiable , which is a common assumption for many NLP optimisation algorithms[70]. For convenience, we introduce some notations that are used in this work. We define an index set of binary variables as $\mathcal{S} = \{s \mid s = 1, 2, \dots, m\}$. The active nodes during B&B are included in a set $\mathcal{L} = \{n_0, n_1, n_2, \dots\}$, where $n_i$ denotes the $i^{\text{th}}$ node that needs to be explored in B&B. Before solving the NLP subproblem at a node $n_i$, we know the values of $m_F^i$ binary variables, whose indices are represented in a set $\mathcal{S}_F^i = \{d_1, d_2, \dots, d_{m_F^i}\}$. These binary variables are $\mathbf{y}_F^i = \left[y_{d_1}, y_{d_2}, \dots, y_{d_{m_F^i}}\right]^T$. The remaining $m_R^i$ binary variables are unknown and relaxed to 0-1 continuous variables that needs to be determined. The indices of these binary variables are included into a set $\mathcal{S}_R^i = \{l_1, l_2, \dots, l_{m_R^i}\}$. These binary variables are $\mathbf{y}_R^i = \left[y_{l_1}, y_{l_2}, \dots, y_{l_{m_R^i}}\right]^T$. Therefore, we have $\mathcal{S}_F^i \cup \mathcal{S}_R^i = \mathcal{S}$ and $m_F^i + m_R^i = m$. We also use $f^{i,*}$ to denote the optimal objective value of the NLP subproblem at node $n_i$, $f^{ub}$ to denote the best upper bound of the MINLP problem until the current iteration. Note that the upper bound is a feasible integer solution for the original MINLP problem P.

**Solution approach**

Compared with other convex MINLP algorithms, the nonlinear B&B algorithm[71] usually has better performance in convergence for solving the strongly nonconvex MINLP problem $\mathbf{P}$[52,68]. This is because the NLP subproblem that needs to be solved at a node $\boldsymbol{n_i}$ differs from that at its parent node in one binary variable only. Threfore, the solution from the parent node is possibly a good initial point for those at the child nodes, while good initial point is rarely availiable during the iterations of other algorithms[68]. The classical nonlinear B&B algorithm for an MINLP problem is presented in **Algorithm S1** of the **Supplementary Material**.

However, when the models are large and strongly nonlinear, the solutions from parent and child nodes may still differ a lot in a nonlinear B&B algorithm, so no good initial point can be provided for the child nodes, causing divergence of the NLP subproblems, as reported in Flores-Tlacuahuac & Biegler for integrated design and control[52] and in Ma et al. for process synthesis[37]. To resolve such problem, we can start from the solution at the parent node and gradually approach the solution of the NLP subproblem at node $\boldsymbol{n_i}$. This strategy can be achieved by using the homotopy continuation (HC) method, which is just a method that starts from the solution of a problem that is easier to solve and then approaches the solution of the original problem gradually[72]. The HC method has been widely used to solve difficult nonlinear algebraic equation systems[73,74] and NLP optimisation problems[66,67] where a homotopy path is constructed through modifying the KKT conditions of the NLP problem. In the sequel, we introduce the HC method in detail and use it to construct a new homotopy path through exploiting the similarity of the NLP subproblems during B&B to enhance convergence in solving the NLP subproblem at a node. Therefore, the new homotopy path is completely different from those in the literature[66,67].

**Homotopy continuation method**

An NLP subproblem to be solved at a node $n_i$ in the nonlinear B&B algorithm is defined below,

$$\min_{\mathbf{x}, \mathbf{y}_R^i} f\left(\mathbf{x}, \mathbf{y}_R^i; \mathbf{y}_F^i\right) \tag{NLP0}$$

$$s.t.\ \mathbf{h}\left(\mathbf{x}, \mathbf{y}_R^i; \mathbf{y}_F^i\right) = 0$$

$$\mathbf{g}\left(\mathbf{x}, \mathbf{y}_R^i; \mathbf{y}_F^i\right) \leq 0$$

$$\mathbf{x} \in \mathbb{R}^n, \mathbf{y}_R^i \in \mathbb{B}_{0,1}^{i,R}, \mathbf{y}_F^i \in \mathbb{Z}^{m_F^i}$$

where $\mathbb{B}_{0,1}^{i,R} := \prod_1^{m_R^i}[0,1] \subset \mathbb{R}^{m_R^i}$ is a rectangular region, and $\mathbf{y}_F^i$ are binary variables with known and fixed values a priori. The NLP subproblem at a node $n_i$ is denoted as $\text{NLP0}(\mathbf{x}, \mathbf{y}_R^i; \mathbf{y}_F^i)$ with the optimal solution notated as $\left(\mathbf{x}^{i,*}, \mathbf{y}_R^{i,*}\right)$. Due to the characterisitcs of the B&B algorithm, the NLP subproblem at node $n_i$ has slight difference in $\mathbf{y}_F^i$ and $\mathbf{y}_R^i$ with that at its corresponding parent node (denoted as $n_{pi}$). Hence, we can use the solution obtained at node $n_{pi}$ to gradually approach the solution at node $n_i$, leading to the construction of homotopy paths. After the NLP subproblem at the parent node $n_{pi}$ is solved, a variable in $\mathbf{y}_R^{pi,*}$ with a fraction value is selected for branching. We use $bpi$ to denote the index of the binary variable with a fractional value for branching at node $n_{pi}$. It means $y_{bpi}^{pi,*}$ is used for branching and hence two child nodes $n_i$ and $n_j$ are created. At the child node $n_i$, the variable in $\mathbf{y}_R^{pi}$ whose index is $bpi$ (i.e., $y_{bpi}^i$) has a known and fixed value at 0 or 1. Therefore, it is removed from $\mathbf{y}_R^i$ but included in $\mathbf{y}_F^i$ when the NLP subproblem at node $n_i$ is solved. Similarly, the binary variable $y_{bpi}^j$ at the child node $n_j$ has a known and fixed value at 1 or 0. Therefore, it is also removed from $\mathbf{y}_R^j$ but included in $\mathbf{y}_F^j$. This is illustrated in Fig. 1.

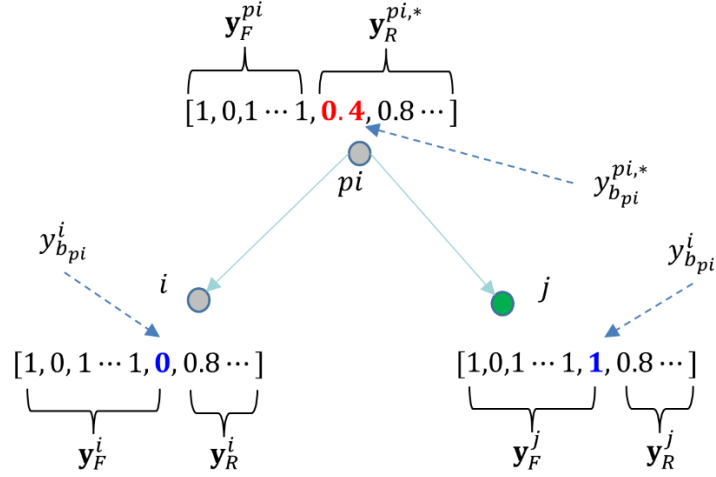

Figure 1 Schematic of branching the parent node $n_{pi}$ into two children nodes $n_i$ and $n_j$.

To ensure feasibility or optimality of the NLP subproblem at node $n_i$, we gradually approach the value of $y_{bpi}^i$ at node $n_i$ (i.e., 0 or 1) from the value of $y_{bpi}^{pi,*}$ using a homotopy parameter $t \in [0,1]$, as shown in Eq. (1).

$$\tilde{y}_{bpi}^i(t) = (1-t) \cdot y_{bpi}^{pi,*} + t \cdot y_{bpi}^i \qquad\qquad t \in [0,1] \qquad\qquad (1)$$

where $\tilde{y}_{bpi}^i(t)$ is called the homotopy variable. When $t = 0$, $\tilde{y}_{bpi}^i(t) = y_{bpi}^{pi,*}$, whilst $\tilde{y}_{bpi}^i(t) = y_{bpi}^i$ when $t = 1$.

Based on this, for a given $t$, the fixed variables $\tilde{\mathbf{y}}_F^i(t)$ with known values at node $n_i$ can be generally represented as follows,

$$\tilde{\mathbf{y}}_F^i(t) = \mathbf{y}_F^i + (1-t)\left(y_{bpi}^{pi,*} - y_{bpi}^i\right)\boldsymbol{e}_{bpi} \qquad\qquad (2)$$

where $\boldsymbol{e}_{bpi}$ is a column vector with the dimension of $m_F^i$ whose entries are all 0 except 1 at $bpi$.

An NLP subproblem with respect to the parameter $t$ ($t \in [0,1]$) can be constructed as follows,

$$\min_{\mathbf{x}, \mathbf{y}_R^i} f\left(\mathbf{x}, \mathbf{y}_R^i; \tilde{\mathbf{y}}_F^i(t)\right) \qquad\qquad \text{(NLPFX)}$$

$$s.t.\ \mathbf{h}\left(\mathbf{x}, \mathbf{y}_R^i; \tilde{\mathbf{y}}_F^i(t)\right) = 0$$

$$\mathbf{g}\left(\mathbf{x}, \mathbf{y}_R^i; \tilde{\mathbf{y}}_F^i(t)\right) \leq 0$$

$$\mathbf{x} \in \mathbb{R}^n, \mathbf{y}_R^i \in \mathbb{B}_{0,1}^{i,R}.$$

For convenience, the parametric optimisation problem is notated as $\text{NLPFX}\big(\mathbf{x}, \mathbf{y}_R^i; t\big)$. Clearly, when $t = 0$ the optimal solution $\big(\mathbf{x}^{pi,*}, \mathbf{y}_R^{pi,*}\big)$ obtained at node $n_{pi}$ is also the optimal solution of the subproblem $\text{NLPFX}\big(\mathbf{x}, \mathbf{y}_R^i; t\big)$. When $t = 1$, the problem $\text{NLPFX}\big(\mathbf{x}, \mathbf{y}_R^i; t\big)$ is equivalent to problem $\text{NLP0}\big(\mathbf{x}, \mathbf{y}_R^i; \mathbf{y}_F^i\big)$. Let $F_{FX}^i(t)$ denote the mapping from $t$ to the optimum of problem $\text{NLPFX}\big(\mathbf{x}, \mathbf{y}_R^i; t\big)$. In other words,

$$F_{FX}^i(t) = \min_{\mathbf{x}, \mathbf{y}_R^i} f\left(\mathbf{x}, \mathbf{y}_R^i; \tilde{\mathbf{y}}_F^i(t)\right)$$

$$s.t. \ \mathbf{h}\big(\mathbf{x}, \mathbf{y}_R^i; \tilde{\mathbf{y}}_F^i(t)\big) = 0$$

$$\mathbf{g}\big(\mathbf{x}, \mathbf{y}_R^i; \tilde{\mathbf{y}}_F^i(t)\big) \leq 0$$

$$\mathbf{x} \in \mathbb{R}^n, \mathbf{y}_R^i \in \mathbb{B}_{0,1}^{i,R}$$

We assume $F_{FX}^i(t)$ is continuous with respect to $t$ within the interval $t \in [0,1]$, which is a mild assumption and satisfied in most cases. Then, the optimal solution of $\text{NLPFX}\big(\mathbf{x}, \mathbf{y}_R^i; t^\nu\big)$ provides a good initial point for problem $\text{NLPFX}\big(\mathbf{x}, \mathbf{y}_R^i; t^{\nu+1}\big)$ if $t^\nu$ and $t^{\nu+1}$ are close enough, i.e., $\varDelta t^\nu := t^{\nu+1} - t^\nu$ is small enough where $\nu$ denotes the iteration number. The homotopy paths of the optimum $F_{FX}^i(t)$ and $F_{FX}^j(t)$ from solving problem $\text{NLPFX}\big(\mathbf{x}, \mathbf{y}_R^i; t\big)$ and $\text{NLPFX}\big(\mathbf{x}, \mathbf{y}_R^j; t\big)$ at nodes $n_i$ and $n_j$ created from the same parent node are illustrated in Fig. 2 (left).

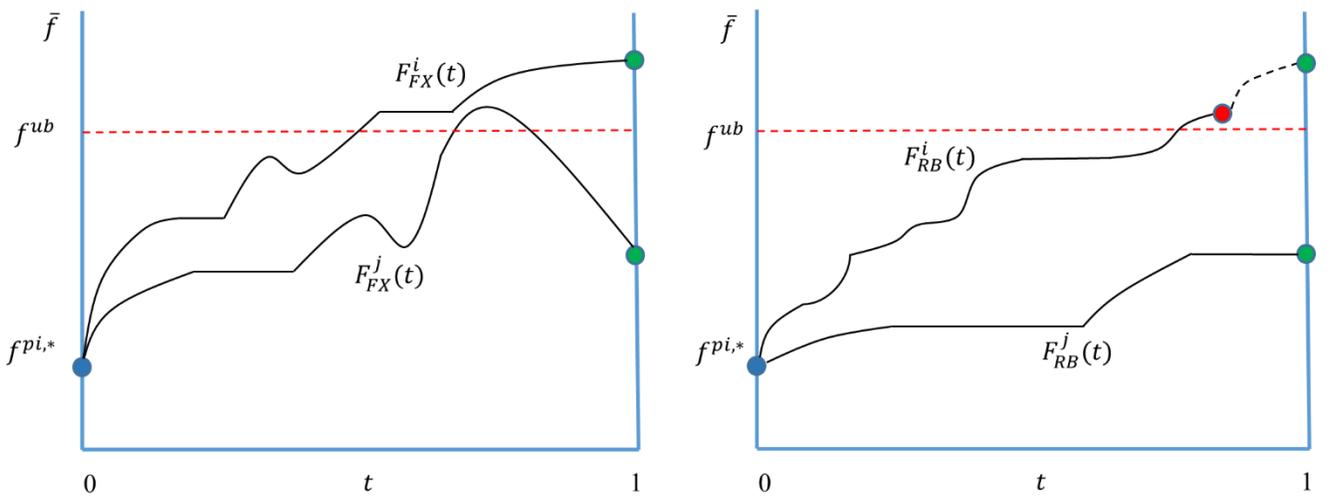

Figure 2 Homotopy paths when solving NLPFX (left) and NLPRB (right) at nodes $n_i$ and $n_j$.

The subproblem $\text{NLPFX}\left(\mathbf{x}, \mathbf{y}_R^i; t\right)$ needs to be solved to optimality from $t = 0$ until $t = 1$, which could significantly increase the computational expense in some cases. As we only need to obtain the optimal solution when $t = 1$ at a node, it is unnecessary to generate the optimal solution for any $t \in (0, 1)$. Based on this insight, a variant of applying the HC method is proposed which first solves a series of NLP feasibility problems (denoted as NLPFP) with a constant objective function (e.g., 0) for $t \in (0, 1]$ during HC and then solve the problem NLP0 to optimality when $t = 1$ at node $n_i$.

$$\underset{\mathbf{x}, \mathbf{y}_R^i}{\text{Min}} \, c \qquad\qquad\qquad \text{(NLPFP)}$$

$$s.t. \, \mathbf{h}\left(\mathbf{x}, \mathbf{y}_R^i; \tilde{\mathbf{y}}_F^i(t)\right) = 0$$

$$\mathbf{g}\left(\mathbf{x}, \mathbf{y}_R^i; \tilde{\mathbf{y}}_F^i(t)\right) \leq 0$$

$$\mathbf{x} \in \mathbb{R}^n, \mathbf{y}_R^i \in \mathbb{B}_{0,1}^{i,R}$$

where $c$ is a constant. The problem NLPFP is specifically notated as $\text{NLPFP}\left(\mathbf{x}, \mathbf{y}_R^i; t\right)$.

When the above two HC variants solve NLP subproblems from $t = 0$ to $t = 1$, both of them cannot use the information of the current lower or upper bounds identified to determine if the current node is fathomed during homotopy. This is because the optimum of both problems NLPFX and NLPFP does not increase monotonically with $t$. A more advantageous variant is designed to gradually tighten the bound of $\tilde{y}_{bpi}^i(t)$ with increase of $t$ so that the optimum monotonically increases with $t$. Once the optimal objective value is greater than the current upper bound (minimization problems) at some $t$ ($t < 1$), it is not needed to continue the HC to approach $t = 1$. In other words, this node is fathomed. The NLP subproblem that needs to be solved in this variant is provided below, which is denoted as NLPRB.

$$\min_{\mathbf{x}, \mathbf{y}_R^i, \tilde{y}_{bpi}^i} f\left(\mathbf{x}, \mathbf{y}_R^i, \tilde{y}_{bpi}^i; \mathbf{y}_{F,s}^i\big|_{s \neq bpi}\right) \qquad \text{(NLPRB)}$$

$$s.t. \quad \mathbf{h}\left(\mathbf{x}, \mathbf{y}_R^i, \tilde{y}_{bpi}^i; \mathbf{y}_{F,s}^i\big|_{s \neq bpi}\right) = 0$$

$$\mathbf{g}\left(\mathbf{x}, \mathbf{y}_R^i, \tilde{y}_{bpi}^i; \mathbf{y}_{F,s}^i\big|_{s \neq bpi}\right) \leq 0$$

$$(1-t) \cdot y_{bpi}^{pi,*} + t \cdot y_{bpi}^i \leq \tilde{y}_{bpi}^i \leq 1, \qquad \text{if } y_{bpi}^i = 1$$

$$0 \leq \tilde{y}_{bpi}^i \leq (1-t) \cdot y_{bpi}^{pi,*} + t \cdot y_{bpi}^i, \qquad \text{if } y_{bpi}^i = 0$$

$$\mathbf{x} \in \mathbb{R}^n, \mathbf{y}_R^i \in \mathbb{B}_{0,1}^{i,R}$$

where $\mathbf{y}_{F,s}^i\big|_{s \neq bpi}$ denotes a vector of the binary variables in $\mathbf{y}_F^i$ at node $n_i$ whose indices are not $bpi$. For convenience, the problem NLPRB for a given $t$ is specifically notated as $\text{NLPRB}\left(\mathbf{x}, \mathbf{y}_R^i, \tilde{y}_{bpi}^i; t\right)$. Similar to problem NLPFX, when $t = 0$, the optimal solution $\left(\mathbf{x}^{pi,*}, \mathbf{y}_R^{pi,*}\right)$ obtained at node $n_{pi}$ is also the optimal solution of the subproblem $\text{NLPRB}(\mathbf{x}, \mathbf{y}_R^i, \tilde{y}_{bpi}^i; t)$. When $t = 1$, problem $\text{NLPRB}\left(\mathbf{x}, \mathbf{y}_R^i, \tilde{y}_{bpi}^i; t\right)$ is equivalent to $\text{NLP0}\left(\mathbf{x}, \mathbf{y}_R^i; \mathbf{y}_F^i\right)$. Let $F_{RB}^i(t)$ denotes the mapping from $t$ to the optimum of problem $\text{NLPRB}\left(\mathbf{x}, \mathbf{y}_R^i, \tilde{y}_{bpi}^i; t\right)$. With $t$ approaching 1, the domain of $\tilde{y}_{bpi}^i$ contracts gradually, leading to the feasible region of the subproblem NLPRB smaller. Hence, $F_{RB}^i(t)$ increases monotonically with $t$, assuming global optimal solution is obtained at each HC step. That is $F_{RB}^i(t^{\nu+1}) \geq F_{RB}^i(t^\nu)$ for any $t^\nu, t^{\nu+1} \in [0,1]$. If $F_{RB}^i(t)$ is larger than the upper bound of the MINLP problem during HC, it is no need to continue, and the current node is pruned. This is because no better solution could be found at the current node. Thus, the advantage of solving problem $\text{NLPRB}\left(\mathbf{x}, \mathbf{y}_R^i, \tilde{y}_{bpi}^i; t\right)$ is that it is possible to terminate the HC calculation earlier. The homotopy paths [i..e, functions $F_{RB}^i(t)$ and $F_{RB}^j(t)$] when solving problems $\text{NLPRB}\left(\mathbf{x}, \mathbf{y}_R^i, \tilde{y}_{bpi}^i; t\right)$ and $\text{NLPRB}\left(\mathbf{x}, \mathbf{y}_R^i, \tilde{y}_{bpj}^j; t\right)$ are shown in Fig. 2 (right).

**Adaptive step lengths during HC**

In the HC method, the step length $\Delta t^\nu$ $(\Delta t^\nu = t^{\nu+1} - t^\nu)$ can significantly affect the performance of the HC method. This is because if $\Delta t^\nu$ is too large, the previous solution at $t^\nu$ may not be a good initial point for the NLP subproblem at $t^{\nu+1}$ during homotopy, which could lead to infeasibility at $t^{\nu+1}$ and finally lead to infeasibility at $t = 1$. On the other hand, if $\Delta t^{\nu}$ is too small, many homotopy steps are required, resulting in increasing computational effort. To balance convergence and computational efficiency, some strategies proposed below are used to update the HC step length $\Delta t$ adaptively.

**Strategy 1**: When a step length leads to feasibility or optimality of the NLP problem at the iteration $\nu$, it indicates that this step length is reliable and continues to be used for the next iteration $(\nu + 1)$. In other words, the step length at the next iteration can be updated using $\Delta t^{\nu+1} \leftarrow \Delta t^{\nu}$.

**Strategy 2**: When two consecutive NLP problems at iterations $(\nu - 1)$ and $\nu$ are solved using the same step length, it indicates a larger step length may also lead to a converged solution of the NLP problem at iteration $(\nu + 1)$. Therefore, we use a larger step length by updating $\Delta t^{\nu+1} \leftarrow 2\,\Delta t^{\nu}$.

**Strategy 3**: If the NLP subproblem at the current iteration $\nu$ fails with a step length, it means the step length is too large and should be reduced, so we set $\Delta t^{\nu+1} \leftarrow \frac{\Delta t^{\nu}}{2}$.

With the above three strategies, the HC algorithm using adaptive step lengths at a node $n_i$ is proposed in the following **Algorithm 1**. In the algorithm, the initial step length $\Delta t^0$ is set as 0.5. $t^0 = 0$ and $t^1 = 0.5$, and we start from $\nu = 1$. Note that $t^1 = 0.5$ because **Algorithm 1** will only be activated when we cannot solve the original problem $\text{NLP0}\big(\mathbf{x}, \mathbf{y}_R^i; \mathbf{y}_F^i\big)$.

In some situations many small step lengths may be required to solve NLP subproblems, which may increase computational effort. We notice that we may have similar information such as the same homotopy variable with the same required value at some nodes in the B&B. It indicates the step lengths at a previous node could be used when a series of NLP subproblems are solved at a latter node with similar information. With this, we may significantly reduce the effort for finding an appropriate step length. For instance, we have the same homotopy variable with the same required value at nodes $n_i$ and $n_j$ $(j < i)$. At node $n_j$ we already identify suitable step lengths, which can be adopted at node $n_i$ and hence reduce the computational effort in finding suitable step lengths at node $n_i$. To identify nodes with similar information, the following conditions are proposed.

**C1:** $\mathcal{J}^i = \{j \mid bpj = bpi \text{ and } y_{bpj}^j = y_{bpi}^i, \; j < i\}$ is nonempty

**C2:** There exits $j \in \mathcal{J}^i$ such that $\Delta y_b^{i,j} \coloneqq \left| y_{bpj}^{pj,*} - y_{bpi}^{pi,*} \right| < \delta$;

**C3:** $j \leftarrow \operatorname{argmin}_{j \in \mathcal{J}^i} \Delta y_b^{i,j}$

---

**Algorithm 1**: HC algorithm using adaptive step length at node $n_i$

---

    **Data**: $0 < \Delta t_{min} < 1, N > 1, \mathbf{x}^{pi,*}, \mathbf{y}_R^{pi,*}, \mathbf{y}_F^{pi}, \mathbf{y}_R^i, \mathbf{y}_F^i, f^{ub}$

1   **initialization**: $\nu \leftarrow 1, t^0 \leftarrow 0, \Delta t^0 \leftarrow 0.5, \; t^1 \leftarrow 0.5$, set $(\mathbf{x}^{pi,*}, \mathbf{y}_R^{pi,*}, \mathbf{y}_F^{pi})$ as the initial
     point;

2   **while** $\nu < N$ and $\Delta t^{\nu-1} > \Delta t_{min}$ **do**

3       solve NLPRB($\mathbf{x}, \mathbf{y}_R^i, \tilde{y}_{bpi}^i; t^\nu$) or NLPFP($\mathbf{x}, \mathbf{y}_R^i; t^\nu$) from the given initial point; the
         optimal solution has $\mathbf{x}^{\nu,*}, \mathbf{y}_R^{\nu,*}$ and $f^{\nu,*}$;

4       **if** a feasible or optimal solution is generated **then**

5           **if** $t = 1$ **then** break;

7           **else if** the solution is optimal and the optimal value $f^{\nu,*} > f^{ub}$ **then** break;

9           **else if** $\nu = 1$ or $\Delta t^{\nu-1} \neq \Delta t^{\nu-2}$ **then** $\Delta t^\nu \leftarrow \Delta t^{\nu-1}$;

11          **else** $\Delta t^\nu \leftarrow 2 \, \Delta t^{\nu-1}$;

13          $t^{\nu+1} \leftarrow \min(t^\nu + \Delta t^\nu, 1)$, the initial point $\leftarrow$ the current solution;

14       **else**

15          $\Delta t^\nu \leftarrow \frac{\Delta t^{\nu-1}}{2}, t^{\nu+1} \leftarrow t^{\nu-1} + \Delta t^\nu$, initial point $\leftarrow$ the last converged solution;

16      $\nu = \nu + 1$;

17   **if** NLPFP($\mathbf{x}, \mathbf{y}_R^i; t^\nu$)  is solved and a feasible solution is obtained **then**

18      solve the original problem NLP0 and get the final solution ($\mathbf{x}^{i,*}, \mathbf{y}_R^{i,*}, f^{i,*}$);

---

C1 ensures that node $n_j$ has the same homotopy variable with the same required value as node $n_i$. C2 enforces the homotopy variable at node $n_j$ starts from a value $y_{bpj}^{pj,*}$ which is different from the starting point of the homotopy variable at node $n_i$ by at most a small constant $\delta$. It means the homotopy variable at nodes $n_i$ and $n_j$ have a close starting point for HC. C3 is used to select a node with the closest starting point to that at node $n_i$ if there are several nodes satisfying C1 and C2. When node $n_j$ satisfies C1-C3, it is expected that the HC step lengths identified at node $n_j$ can also successfully solve a series of NLP subproblems at node $n_i$ during homotopy, so those steps will be used during the homotopy calculation at node $n_i$. However, once a step length leads to infeasibility of an NLP subproblem, it is then halved and the step length updating follows the same strategies as those in Algorithm 1. The HC algorithm using adaptive step lengths and C1-C3 at node $n_i$ is proposed in **Algorithm 2** as follows,

---

**Algorithm 2**: HC algorithm using variable step length and C1-C3 at a node $n_i$

**Data**: $0 < \Delta t_{min} < 1, N > 1, \mathbf{x}^{pi,*}, \mathbf{y}_R^{pi,*}, \mathbf{y}_F^{pi}, \mathbf{y}_R^i, \mathbf{y}_F^i, \mathcal{T}$;

1  **initialization**: $\nu \leftarrow 1$, set $(\mathbf{x}^{pi,*}, \mathbf{y}_R^{pi,*}, \mathbf{y}_F^{pi})$ as the initial point;

2  **if** there exists $j$ satisfying conditions (1)-(3) **then**

3      $t^1 \leftarrow t^{j,1}$ ($t^1$ at node $j$);

4      **while** True **do**

5          solve NLPRB($\mathbf{x}, \mathbf{y}_R^i, \tilde{y}_{bpi}^i; t^\nu$) or NLPFP($\mathbf{x}, \mathbf{y}_R^i; t^\nu$) from the given initial point; the optimal solution has $\mathbf{x}^{\nu,*}, \mathbf{y}_R^{\nu,*}$ and $f^{\nu,*}$;

6          **if** the solution is feasible or optimal **then**

7              **if** $t^\nu = 1$ **then** break;

7              **else if** the solution is optimal and the optimum $f^{\nu,*} > f^{ub}$ **then** break;

9              **else** $t^{\nu+1} \leftarrow t^{j,\nu+1}$ ($t^{\nu+1}$ at node $j$) , initial point $\leftarrow$ the current solution;

11          **else**

12              $\Delta t^\nu \leftarrow \frac{\Delta t^{j,\nu-1}}{2}, t^{\nu+1} \leftarrow t^{\nu-1} + \Delta t^\nu$, initial point $\leftarrow$ the last converged solution;

13              switch to **Algorithm 1** and start from its line 2; break;

14          $\nu \leftarrow \nu + 1$;

15  **else**

16      apply **Algorithm 1**;

---

**Homotopy continuation enhanced branch and bound algorithm**

The homotopy continuation enhanced branch and bound (HCBB) algorithm is shown in Fig. 3. In the HCBB algorithm, $n_i$ is used to denote node $i$ and $\mathcal{L}$ is the queue of nodes that need to be investigated. In the beginning, queue $\mathcal{L}$ only contains the root node $n_0$ , the upper bound ($f^{ub}$) is initialized as $+\infty$, and all binary variables are included in the set $S_R^i$ with unknown values. At a node $n_i$, **Algorithm 2** is first executed to generate an optimal solution (denoted as $\mathbf{x}^*$, $\mathbf{y}_R^{i,*}$, $f^{i,*}$) or identify infeasibility. If infeasibility is returned, then the latest optimal or feasible solution with the values of the corresponding HC parameter $t$, and the last value of $\Delta t$ during HC are recorded, which will be used for post check and refinement procedure explained in the next subsection. If an optimal solution is identified when solving $\mathrm{NLPRB}\big(\mathbf{x}, \mathbf{y}_R^i, \tilde{y}_{bpi}^i; t^\nu\big)$, then we need to check if the incumbent solution is greater than $f^{ub}$. If greater, then the current node can be discarded as no better optimal solution can be found under this node. Otherwise, we need to check the values of $\mathbf{y}_R^i$ in the solution (i.e., $\mathbf{y}_R^{i,*}$). If all $\mathbf{y}_R^{i,*}$ are 0 or 1, it means we find a feasible integer solution at this node. This node is fathomed and the upper bound is updated if the feasible integer solution is less than $f^{ub}$. Otherwise, branching on this node is conducted. A binary variable in $\mathbf{y}_R^i$ having a fraction value closest to 0.5 is selected to branch, which is a common-used rule[75]. The well-known best-first strategy[76] is used to select a node to be explored in the next iteration, which usually has better performance than the breadth-first strategy[77] and the depth-first strategy[77].

In addition, we allow the solution obtained from a parent node is used as an initial point for the NLP subproblem at its immediate children nodes. In other words, a warm start is used.

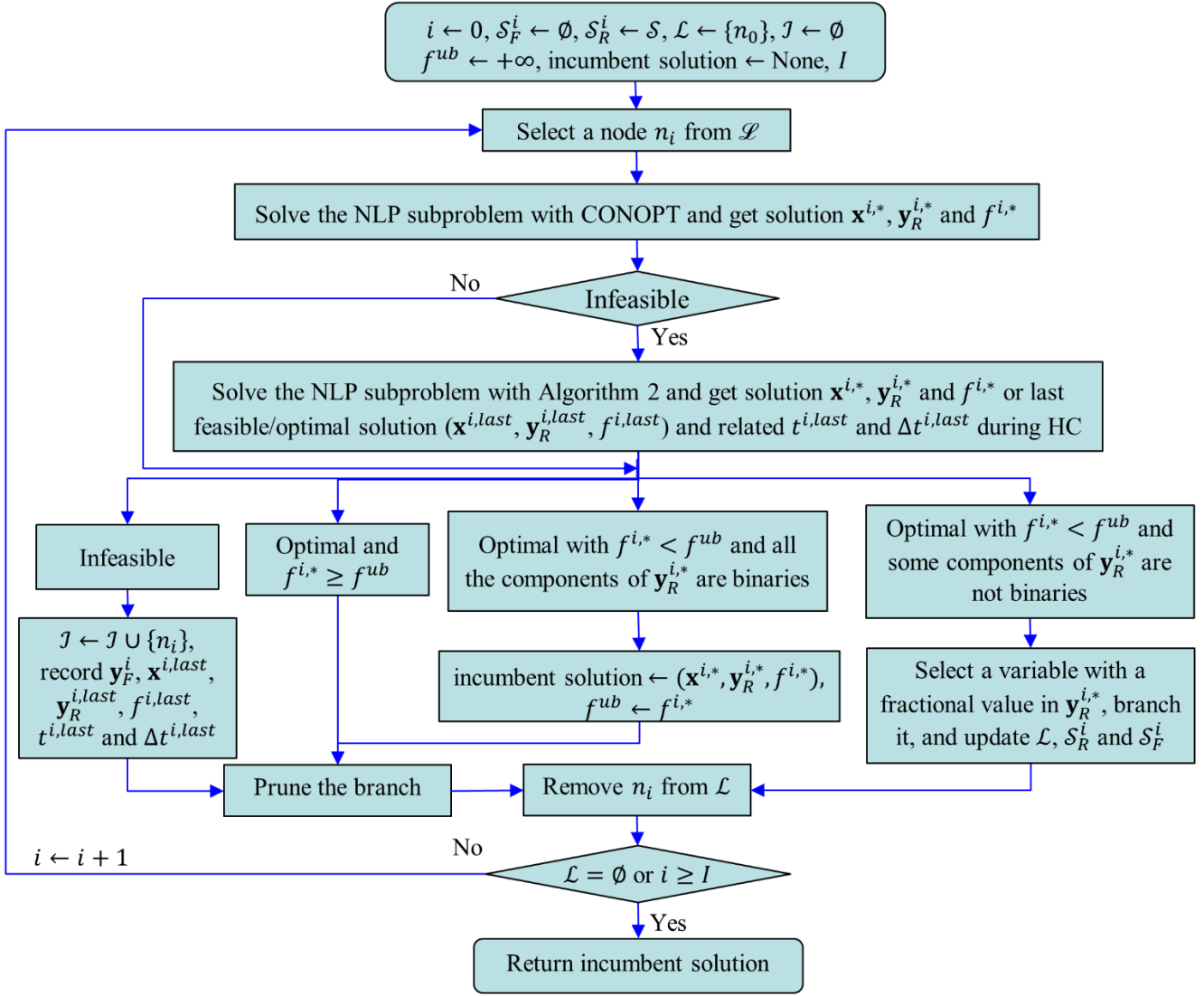

Figure 3 Homotopy continuation enhanced branch and bound algorithm.

**Post check and refinement procedure**

After the MINLP problem is solved using the HC enhanced B&B algorithm in the previous section, the optimal solution (denoted as $\mathbf{x}^*, \mathbf{y}^*, f^*$) may not be the global optimum. One of the reasons is $\Delta t_{\min}$ and the maximum number of HC steps (i.e., $N$) specified are reached at some nodes, leading to no feasible solution found at these nodes, whereas some feasible solution might be found if a larger $N$ and a smaller $\Delta t_{\min}$ are used at these nodes, which may lead to a better solution obtained finally. Based on this consideration, a post check and refinement procedure is conducted at the nodes flagged as infeasibility. We use a set $\mathcal{I}$ to denote those nodes and some corresponding information at each node in $\mathcal{I}$ are also recorded during the above HCBB, which includes $\mathbf{y}_F^i$, the last optimal or feasible solution $(\mathbf{x}_R^{i,last}, \mathbf{y}_R^{i,last}, f^{i,last})$ during the homotopy continuation for an infeasible node $i$, the corresponding HC parameter $(t^{i,last})$, and the last HC step length $(\Delta t^{i,last})$.

---

**Post check and refinement procedure**

---

**Data**: $\mathcal{I}, f^*, N^{post}, \Delta t_{min}^{post}$;

1   **initialization**: $\mathcal{L} \leftarrow \emptyset$;

2   **while** $\mathcal{I} \neq \emptyset$ **do**

3       select a node $n_i$ from $\mathcal{I}$ and let $\mathcal{L} \leftarrow \{n_i\}$, set $f^{ub} = f^*$;

4       **If** NLPRB is solved and $f^{i,last} > f^{ub}$ **then**

5          prune the current node;

6       **else**

7          $N \leftarrow N^{post}, \Delta t_{min} \leftarrow \Delta t_{min}^{post}, t^0 \leftarrow t^{i,last}, \Delta t^0 \leftarrow \Delta t^{i,last}, t^1 \leftarrow t^{i,last} + \Delta t^0$, initial point $\leftarrow (\mathbf{x}^{i,last}, \mathbf{y}_R^{i,last})$;

8          solve the problem NLPRB($\mathbf{x}, \mathbf{y}_R^i, \tilde{y}_{bpi}^i; t$) or NLPFP($\mathbf{x}, \mathbf{y}_R^i; t$) using **Algorithm 2**. The optimal solution has $\mathbf{x}^{i,*}, \mathbf{y}_R^{i,*}$ and $f^{i,*}$;

9          **if** the solution is infeasible **then**

10            prune the current branch;

11          **else**

12            **if** $f^{i,*} \geq f^{ub}$ **then**

13               prune the current branch;

14            **else if** all components of $\mathbf{y}_R^{i,*}$ are binaries **then**

15               incumbent solution $\leftarrow (\mathbf{x}^{i,*}, \mathbf{y}_R^{i,*}, f^{i,*}), f^{ub} \leftarrow f^{i,*}$, prune the branch;

16            **else**

17               select a variable with a fractional value in $\mathbf{y}_R^{i,*}$; generate two nodes through fixing the variable at 0 and 1 respectively; put the two nodes in $\mathcal{L}$;

18               set $N, \Delta t_{min}, t^0, t^1$ and $\Delta t^0$ to default values and set the current solution to be the initial point;

19               conduct B&B again with $\mathcal{I}$ adjusted accordingly;

---

During post check and refinement of HCBB algorithm solving NLPRB subproblems, we first compare the recorded optimum ($f^{i,last}$) at node $i$ in $\mathcal{I}$ with $f^{ub}$ (equal to the current optimum of the MINLP). If it is greater than $f^{ub}$, then this node is discarded and removed from $\mathcal{I}$ directly without solving any new NLP subprobem. This is called a post check step. Otherwise, a larger maximum number of HC steps (denoted as $N^{post}$, e.g., 1000) and a smaller minimum step length (denoted as $\Delta t^{post}_{min}$, e.g., $1 \times 10^{-15}$ which is around 10 times the machine accuracy for double-precision float point arithmetic[78]) are used to solve problem NLPRB starting from the recorded optimal or feasible solution. This is called a refinement step. If no better optimal solution is found for $t = 1$ at this node or the node is still flagged as infeasible due to reaching $N$ or $\Delta t_{min}$, then this node is pruned from $\mathcal{I}$. Otherwise, this node is branched, and two children nodes are created and put into $\mathcal{L}$. The HCBB algorithm in Fig. 3 is used to investigate the new created children nodes again. It should be noted that when NLPFP is solved in HCBB, the refinement step is applied directly because the post check step is not applicable.

The complete HCBB algorithm incorporating post check and refinement procedure is provided in **Appendix A of** the **Supplementary Material**. In the HCBB algorithm, we can solve several NLP variants including NLPFX, NLPFP, and NLPRB during HC at a node $n_i$. According to our previous analysis, solving both NLPFP and NLPRB are more advantageous than solving NLPFX. Thus, we only implement HCBB algorithms with solving NLPFP and NLPRB, which are denoted as HCBB-FP and HCBB-RB respectively. All these HCBB algorithms are implemented in Python[79], which exchanges data with GAMS via Python API of GAMS[80]. Each NLP subproblem is solved using GAMS/CONOPT[81], as it is well-known for good performance.

**Remarks:** 1) The complete HCBB-FP algorithm can theoretically guarantee global optimality if the NLP subproblem for $t = 1$ at each node can be solved to global optimality.

2) The complete HCBB-RB algorithm converges to a local optimum if the NLP subproblems at some $t$ during HC cannot be solved to global optimality.

3) The HCBB algorithm may fail if a feasible or locally optimal solution at root node is not obtained.

**Computational studies**

Four process synthesis problems using rigorous unit operation models are solved to illustrate the capability of the proposed HCBB-FP and HCBB-RB. GAMS/BARON, GAMS/ANTIGONE, GAMS/SBB[68], GAMS/DICOPT[69], SRMINLP from our previous work[37] and the nonlinear B&B implemented in Python by us (denoted as BB) are also used to solve these examples for comparison. Note that our implementation of the nonlinear B&B algorithm may be different from GAMS/SBB[68]. All examples are modelled in GAMS 24.3.3 on a desktop with 3.20 GHz Intel Core i5-3470 CPU, 8GB RAM and 64-bit operating system. As expected, BARON and ANTIGONE cannot generate a feasible solution for any example within 1 hour.

Two initialization strategies are used to generate several initial points for the proposed HCBB algorithms, GAMS/SBB, GAMS/DICOPT, SRMINLP, and BB. These two initialization strategies are similar to that of Ma et al.[37,82]. We decouple the reactor network from the distillation system and solve them separately. We first solve the reactor network, which provides an input for the distillation system. The distillation system is then simulated or optimized to obtain an initial point for the entire process. While the distillation system is simulated without guaranteeing feasibility of all constraints such as product purity constraints in the first initialization strategy, the second strategy applies the hybrid steady-state and time-relaxation-based optimisation algorithm[83] to get a feasible solution of the distillation system. More details are provided in the **Supplementary Material**.

To be clear, the "rigorous models" in these case studies refer to the MESH model for distillation columns and CSTR and PFR models for reactors, as they are frequently used in process simulators. Especially, the bulk of interest in rigorous models is often related to distillation columns[58].

**Example 1: benzene chlorination process**

This example is adopted from Zhang et al.[54] It uses benzene and chlorine to produce desirable product chlorobenzene. The reactions are provided in Fig. S1. The superstructure is presented in Fig. S2. Fresh and recycled raw materials (i.e., benzene and chlorine) are mixed and fed into the reactor network to produce chlorobenzene and byproduct dichlorobenze. The highly volatile hydrochloric acid and chlorine are recovered from the reaction effluent using a flash unit and recycled back to the reaction system. The benzene, chlorobenzene and dichlorobenze mixture from the flash bottom is then heated to the bubble point and fed to the distillation system to obtain pure products chlorobenzene and dichlorobenzene and recover unreacted benzene that is recycled back to the reaction system. The productivity of chlorobenzene should be at least 50 kmol h$^{-1}$. The recovery and purity of chlorobenzene and dichlorobenzene are at least 0.99, and the purity of the recycled benzene should also be at least 0.99.

Table 1 Optimisation results of different algorithms for Example 1

| Item | DICOPT | SRMINLP | SBB | BB | HCBB-FP | HCBB-RB |
|---|---|---|---|---|---|---|
| TAC (M€ year$^{-1}$) | 0.614/Inf /0.614 | 0.614/Inf /0.614 | 0.614/0.614 /0.614 | 0.613/0.614 /0.614 | 0.613/0.614 /0.614 | 0.613/0.614 /0.614 |
| CPU Time (s) | 82/28/71 | 4/6/3 | 5/7/5 | 4/3/5 | 14/14/12 | 7/12/8 |
| $N_{node}$ | 17/17/17 | 2/51/2 | 4/3/3 | 3/3/3 | 3/3/3 | 3/3/3 |
| $N_{inf}$ | 15/16/14 | 0/44/0 | 1/1/1 | 1/1/1 | 1/1/1 | 0/1/0 |
| $N_{nlp}$ | 17/17/17 | 2/51/2 | 4/3/3 | 3/3/3 | 13/13/12 | 5/13/4 |
| $t_{post}$ (s) | - | - | - | - | 215/155/70 | 0/0/0 |
| $N_{inf,post}$ | 15/16/14 | 0/44/0 | 1/1/1 | 1/1/1 | 1/1/1 | 0/0/0 |

Inf: infeasible solution; $N_{node}$: total number of nodes fathomed; $N_{inf}$: number of nodes with infeasibility; $N_{nlp}$: total number of NLP problems.
$t_{post}$: computational time for post check and refinement; $N_{inf,post}$: number of nodes with infeasibility after post check and refinement.

The superstructure is modelled using the MINLP formulation from Ma et al.[37] with minimization of TAC. We use HCBB-FP, HCBB-RB, BB, SBB, DICOPT and SRMINLP to solve the MINLP formulation from three initial points generated. The optimisation results are provided in Table 1 where the results from different initial points are separated using slashes. There are totally 8 binary variables, 2450 continuous variables and 2451 constraints involved. As seen from Table 1, HCBB-FP and HCBB-RB generate almost the same optimal solution with TAC around 0.613 or 0.614 M€ year$^{-1}$ from all three different initial points, which is 1.6% less than that (0.624 M€ year$^{-1}$) from Zhang et al.[54] The B&B algorithms including SBB, BB, HCBB-FP and HCBB-RB have better convergence than DICOPT and SRMINLP as both DICOPT and SRMINLP fail to find a feasible solution from the second initial point, whilst the B&B algorithms can find the optimal solution from all three initial points. DICOPT fails to find a feasible solution from the second initial point because only the NLP problem at the root node ($N_{\text{node}} - N_{\text{inf}} = 1$) is solved to feasibility. SRMINLP also fails from the second initial point because the complementary constraints prohibit the problem from converging to a feasible point. From Table 1, we can also observe that SBB and BB require less computational effort than HCBB-FP and HCBB-RB because the formers only solve one NLP subproblem at a node, while the latters need to solve several NLP subproblems when directly solving the original NLP subproblem fails, as indicated by the number of NLP subproblems solved ($N_{\text{nlp}}$).

From Table 1 it can also be observed that HCBB-FP needs to solve more NLP subproblems than HCBB-RB from the first and third initial points (13 vs. 5 and 12 vs. 4). This attributes to the fact that once a solution from HCBB-RB during HC at a node has an objective value greater than the upper bound, this node is flagged as having an optimal value larger than the upper bound, whilst such strategy cannot be implemented in HCBB-FP. For example, HCBB-FP and HCBB-RB find an incumbent solution with TAC of 0.614 M€ year$^{-1}$ at node 1 from the initial point 1. At node 2, HCBB-FP solves the NLP subproblems until the step length $\Delta t = 0.008$ is less than the minimum allowable step length (i.e., 0.01) during HC. However, HCBB-RB obtains a solution with TAC of 0.687 M€ year$^{-1}$ at $t = 0.25$, which is greater than the upper bound (i.e., 0.614 M€ year$^{-1}$), leading to termination of the HC calculation. Similar observation can be made for HCBB-FP and HCBB-RB from the initial point 3. However, the situation is different from the initial point 2. The computational time for HCBB-FP and HCBB-RB is very close (14 vs. 12). This is because there is no incumbent solution found yet when the node 1 is being solved and hence HCBB-RB has no chance to terminate HC early. Both algorithms cannot find an optimal solution at node 1 even if the step length less than $\Delta t_{min}$ specified. As a result, node 1 is marked as infeasible in HCBB-FP and HCBB-RB and requires to be investigated again in the post check and refinement, where HCBB-RB finds that the last optimal solution during HC for node 1 with TAC of 0.690 M€ year$^{-1}$ is greater than the optimal value of 0.613 M€ year$^{-1}$. Therefore, there is no need for HCBB-RB to start the refinement. On the contrary, HCBB-FP has to conduct refinement directly, leading to high computational effort. As seen from Table 1, HCBB-FP spent additional 155 CPU seconds to investigate the infeasible node 1. From this example, we can conclude that HCBB-RB is better than HCBB-FP, whilst SBB and BB are superior to HCBB-RB.

In the optimal design generated from HCBB-RB, a few bypass efficiencies have fractional values probably due to numerical errors or local optimum obtained. These fractional values are rounded up and an NLP problem with all binary variables and bypass efficiencies fixed at 0 or 1 is then solved. The final TAC changes little relative to that from HCBB-RB, (i.e., not more than 0.01%). The final optimal design is shown in Fig. S3.

**Example 2: cyclohexane oxidation process**

This example is from Zhang et al.[54] Cyclohexane is used to produce a mixture of cyclohexanol and cyclohexanone (i.e., KA oil) through oxidation. The reaction pathways are provided in Fig. S4 where the desired reaction produces cyclohexanol and cyclohexanone and the side reaction produces adipic acid. The superstructure is illustrated in Fig. S5. Air and cyclohexane are fed to the reactor network with only CSTR reactors used. The liquid effluent drawn from the reactor network is separated using two distillation columns to generate the desired product and byproduct. The recovered cyclohexane from the distillation system is recycled to the reactor network. The purity and recovery of products and raw material should be at least 0.9995.

The superstructure is modelled using the MINLP formulation from Ma et al.[37] in which the Wilson equation is used to calculate liquid activity coefficients, whilst the vapor is assumed to be ideal gas. The objective is to minimize TAC including energy cost and annualized capital cost. The optimisation problem involves 12 binary variables, 4647 continuous variables, and 4646 constraints. We generate three different initial points to initialize all these algorithms. The optimisation results are provided in Table 2 where the results from different initial points are separated using slashes.

Table 2 Optimisation results of different algorithms for Example 2

| Item | DICOPT | SRMINLP | SBB | BB | HCBB-FP | HCBB-RB |
|---|---|---|---|---|---|---|
| TAC (M€ year$^{-1}$) | Inf/Inf /Inf | 2.146/2.118 /2.118 | Inf/Inf /Inf | Inf/Inf /3.572 | 2.146/2.146 /2.146 | 2.118/2.118 /2.118 |
| CPU Time (s) | 24/36/3600+ | 13/24/11 | 25/24/15 | 59/227/48 | 187/554/105 | 115/377/109 |
| $N_{node}$ | 7/20/8 | 2/2/2 | 3/3/3 | 15/62/11 | 15/55/13 | 15/55/15 |
| $N_{inf}$ | 6/19/7 | 0/1/0 | 2/2/2 | 8/28/4 | 4/15/0 | 0/0/1 |
| $N_{nlp}$ | 7/20/8 | 2/2/2 | 3/3/3 | 15/62/11 | 89/289/46 | 37/144/36 |
| $t_{post}$ (s) | - | - | - | - | 3600+/3600+/0 | 0/0/0 |
| $N_{inf,post}$ | 6/19/7 | 0/1/0 | 2/2/2 | 8/28/4 | 4/15/0 | 0/0/0 |

Inf: infeasible solution; $N_{node}$: total number of nodes; $N_{inf}$, number of nodes with infeasibility. $N_{nlp}$: total number of NLP problems.
$t_{post}$: computational time for post check and refinement; $N_{inf,post}$: number of nodes with infeasibility after post check and refinement.

As can be seen from Table 2, neither DICOPT nor SBB can solve the problem from any initial point because the NLP subproblems become more difficult to solve when the nonideal physical property Wilson equation is used. This can be evidenced from the fact that only the NLP subproblem at the root node is solved in both DICOPT and SBB. BB can only obtain a worse local optimum from the third initial point, whilst it fails from the other two initial points. HCBB-FP and HCBB-RB can get the optimal solution from all three initial points, showing their robustness in convergence. This indicates the effectiveness of the HC method in the improvement of convergence. More importantly, HCBB-RB converges to a better local optimum of 2.118 M€ year$^{-1}$, which is 1% less than the locally optimal solution of 2.146 M€ year$^{-1}$ obtained by HCBB-FP. SRMINLP solves the problem to local optimality of 2.118

M€ year$^{-1}$ from the second and third initial points and 2.146 M€ year$^{-1}$ from the first initial point , which is worse than that from HCBB-RB. However, SRMINLP requires less computational time than HCBB-FP and HCBB-RB due to far fewer NLP subproblems being solved.

We now compare HCBB-FP and HCBB-RB in detail. From Table 2, we can observe HCBB-RB uses 26-52% less computational time than HCBB-FP to locate a locally optimal solution from an initial point. This is because HCBB-FP requires many HC steps to solve NLP subproblems, whilst HCBB-RB can terminate earlier once a larger TAC than the upper bound is obtained. As a result, HCBB-RB solves 45-60% fewer NLP subproblems than HCBB-FP. Furthermore, HCBB-RB claimed one infeasible node only from the initial point 3, while HCBB-FP claimed infeasibilities at 4, 14 and 13 nodes from the three initial points, respectively. More importantly, HCBB-FP spends more than 1 hour to complete the post check and refinement, whilst HCBB-RB requires tiny computational time to complete the post check and refinement. This is because HCBB-RB finds the last optimal solution during HC at the infeasible node is larger than the upper bound and hence the refinement procedure is not required.

There are several fractional bypass efficiencies in the locally optimal solution. We round them up and solve the resulted NLP problem again to obtain a slightly higher TAC of 2.120 M€ year$^{-1}$, which is still 4% − 5% less than the optimal TAC from Zhang et al.[54] (2.23 M€ year$^{-1}$). The optimal design is shown in Fig. S6. The lower TAC is because we invest a little more in reactors to achieve higher conversion and selectivity while decreasing the burden of separation systems evidently.

**Example 3: hydrodealkylation process of toluene (HDA) using lumped reactor model**

This example is taken from Ma et al[37], which uses toluene and hydrogen to produce product benzene and byproduct diphenyl. The reactions are given in Fig. S7. The superstructure is illustrated in Fig. S8. All data are given in the **Supplementary Material**. The desired benzene molar purity is 99.97% with a production rate of 124.8 kmol h$^{-1}$. The objective is to maximize the economic profit computed by the revenue from benzene and diphenyl minus annualized capital cost and operating cost.

The superstructure is modelled using the MINLP formulation from Ma et al.[37] which is solved using HCBB-FP, HCBB-RB, GAMS/DICOPT, GAMS/SBB, BB, and SRMINLP, respectively. Six different initial points are generated to initialize these algorithms. There are 8142 constraints, 8643 continuous variables, and 13 binary variables in the optimisation problem. The optimisation results are provided in Table 3 where the results from different initial points are separated using slashes. The results from GAMS/DICOPT are not shown in Table 3 due to its infeasibility from any initial point.

Table 3 Optimisation results of different algorithms for Example 3

| Item | SRMINLP | SBB | BB | HCBB-FP | HCBB-RB |
|---|---|---|---|---|---|
| Profit (M$ year⁻¹) | 3.920/4.189 /4.179/4.184 /4.189/4.954 | 4.713/4.710 /4.959/4.958 /4.712/Inf | 4.956/4.954 /4.956/4.958 /4.955/4.958 | 4.956/4.954 /4.956/4.958/ 4.955/4.958 | 4.956/4.954 /4.956/4.958/ 4.955/4.958 |
| CPU Time (s) | 36/41/64 /62/58/18 | 48/66/81 /98/66/25 | 54/88/116 /76/85/63 | 150/170/187 /112/174/123 | 130/148/231 /107/176/220 |
| $N_{node}$ | 2/2/0/2/2/2 | 5/7/7/5/4/5 | 7/5/7/5/7/5 | 9/5/7/5/7/5 | 9/5/9/5/7/5 |
| $N_{inf}$ | 0/0/0/0/0/0 | 2/3/3/2/1/3 | 2/1/2/1/2/1 | 2/1/2/1/2/1 | 0/1/0/1/1/1 |
| $N_{nlp}$ | 2/2/2/2/2/2 | 5/7/7/5/4/5 | 13/10/12 /7/12/22 | 47/23/34 /20/36/33 | 79/24/24 /18/30/88 |
| $t_{post}$ (s) | - | - | - | 1500/656/3600+ /2296/1999/1593 | 0/0/0/0/0/0 |
| $N_{inf,post}$ | 0/0/0/0/0/0 | 2/3/3/2/1/3 | 2/1/2/1/2/1 | 2/1/2/1/2/1 | 0/0/0/0/0/0 |

Inf: infeasible solution; $N_{node}$: total number of nodes; $N_{inf}$: number of nodes with infeasibility; $N_{nlp}$: total number of NLP problems.

$t_{post}$: computational time for post check and refinement; $N_{inf,post}$: number of nodes with infeasibility after post check and refinement.

From Table 3, we can observe that HCBB-FP and HCBB-RB can find the locally optimal solution from any initial point, demonstrating their robustness in convergence. Almost the same locally optimal solution is generated from the six initial points. BB has very similar convergence performance to HCBB-FP and HCBB-RB for this example. SRMINLP and SBB can identify the locally optimal solutions from four out of the six initial points. While SRMINLP fails from the third initial point due to infeasibility caused at the root node, SBB fails from the sixth initial point. Although SRMINLP requires much less computational time than HCBB-FP and HCBB-RB, the optimal profit it generates from some initial points is much less than that obtained from HCBB-FP and HCBB-RB. For instance, SRMINLP generates a locally optimal profit of 3.920 M$ year$^{-1}$ from the first initial point, which is 21% lower than that of 4.956 M$ year$^{-1}$ obtained from HCBB-FP and HCBB-RB from the same initial point. SBB also requires much less computational time than HCBB-FP and HCBB-RB. However, the profit it obtained from some initial points is still much lower than that generated from HCBB-FP and HCBB-RB.

From Table 3 HCBB-FP and HCBB-RB find almost the same local optimum of 4.956 M$ year$^{-1}$ from the six initial points. While HCBB-RB consumes 4%-13% less computational time than HCBB-FP from the first, second and fourth initial points, it requires 1%-79% more computational time from the third, fifth and sixth initial points. HCBB-RB requires to solve much more NLP subproblems from the sixth initial points compared to HCBB-FP. The possible reason is due to the difficulties in solving the problem with stream flow rates appearing in the denominator to optimality when the flow rates approach zero[22,84]. After solving, HCBB-FP and HCBB-RB flag 1-2 nodes as infeasible due to $\Delta t_{\min}$ or the maximum number of HC steps (i.e., $N_{max}$) specified being reached. As a result, HCBB-FP spends 656-3600 CPU seconds in the post check and refinement procedure, whilst HCBB-RB quickly identifies no better solution could be found at these infeasible nodes after the post check with negligible computational effort and hence the time-consuming refinement procedure is not conducted. Therefore, HCBB-RB requires much less total computational effort for B&B and post check and refinement to identify the same optimal solution compared to HCBB-FP.

In the optimal solution directly from HCBB-RB, there are a few fractional bypass efficiencies in the first column. After rounding the bypass efficiencies and optimizing the NLP problems again, the profit becomes 4.956 M$ year$^{-1}$ which changes by around 0.04%, as shown in the **Supplementary Material**. The final optimal design is shown in Fig. 4.

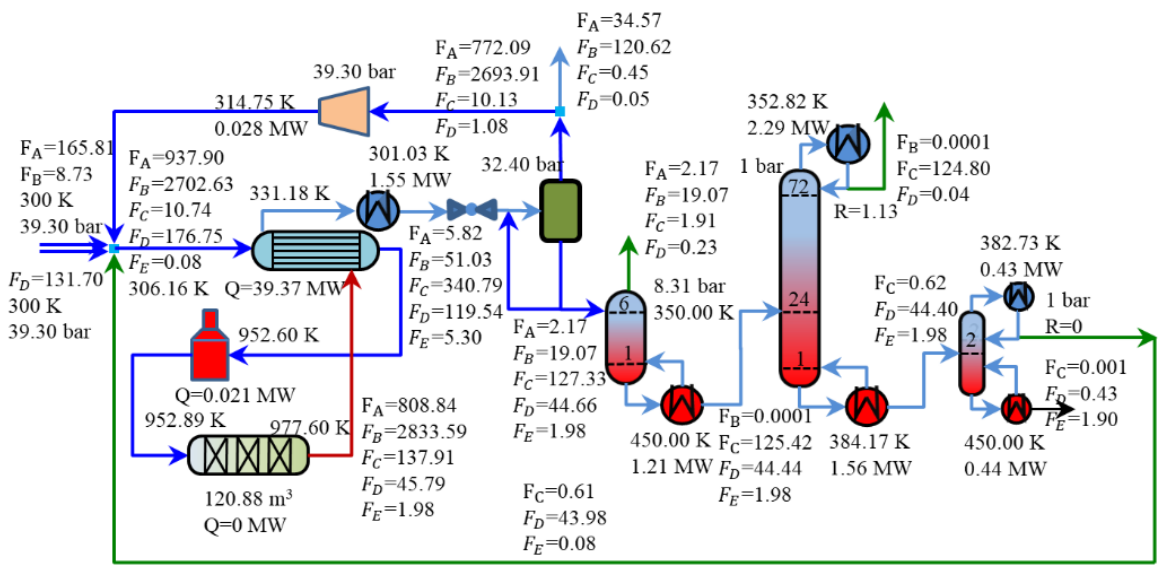

Figure 4 Optimal design of HDA process where the unit of flow rate is kmol h$^{-1}$ and A, B, C, D and E represent $H_2$, $CH_4$, $C_6H_6$, $C_7H_8$ and $C_{12}H_{10}$ respectively.

**Example 4: hydrodealkylation of toluene with differential reactor model**

This example is very similar to Example 3 but the reaction kinetic equation from Dimian et al.[85] is used. Both isothermal and adiabatic reactors are modelled using differential equations and discretized with orthogonal collocation finite element method[86]. The reactor model is more accurate due to the error caused by using an average temperature in the reaction kinetics for the whole adiabatic reactor in Example 3 being avoided. The kinetic equations and differential reactor models after discretization are provided in the **Supplementary Material**, which uses 50 elements and 3 Radau collocation points. Three reactor units are incorporated in the superstructure as shown in Fig. S9.

The optimisation problem involves 24 758 equations, 25 552 continuous variables and 19 binary variables. We generate six different initial points to initialize the algorithms. The computational results are provided in Table 4. As GAMS/DICOPT still cannot solve the problem from any initial point within 1 hour, the results from it are not shown in Table 4. As can be seen from Table 4, SRMINLP often finds a worse local optimum compared to BB, HCBB-FP, and HCBB-RB. For instance, SRMINLP finds a local optimum of 4.291 M\$ year$^{-1}$ from the third initial point, which is 12% lower than the optimal solution of 4.904 M$ year$^{-1}$. SBB finds a locally optimal solution with TAC of 4.291 M$ year$^{-1}$ only from the third initial point, which is much lower than that of 4.904 M$ year$^{-1}$ from HCBB-FP and HCBB-RB. While SBB fails to solve the problem from the other five initial points due to the failure in solving NLP subproblems during B&B, HCBB-FP and HCBB-RB find the optimal solution with TAC of 4.904 M$ year$^{-1}$ from all six initial points. Although BB can find similar optimum as the HCBB algorithms from five initial points, it obtains a much worse local optimum of 4.309 M$ year$^{-1}$ from the sixth initial point compared to that from HCBB-FP and HCBB-RB. It is expected that SRMINLP, SBB and BB require less computational time than HCBB-FP and HCBB-RB due to a smaller number of NLP subproblems being solved.

Table 4 Computational results of different algorithms for Example 4

| Item | SRMINLP | SBB | BB | HCBB-FP | HCBB-RB |
|---|---|---|---|---|---|
| Profit (M$ year$^{-1}$) | 4.309/4.303 /4.291/4.903 /4.639/4.902 | Inf/Inf /4.291/Inf/ Inf/Inf | 4.904/4.904 /4.904/4.904 /4.901/4.309 | 4.904/4.904 /4.904/4.904 /4.901/4.903 | 4.904/4.904 /4.904/4.904 /4.901/4.904 |
| CPU Time (s) | 182/1034/294 /278/198/361 | 132/278/215/ 213/263/263 | 309/475/280/ 269/141/494 | 863/1238/474/ 573/372/1757 | 683/1136/693/ 590/299/1507 |
| $N_{node}$ | 2/2/2/2/2/2 | 3/3/7/7/7/3 | 11/11/7/7/5/11 | 11/13/7/7/5/11 | 11/13/7/7/5/11 |
| $N_{inf}$ | 0/0/0/0/0/0 | 2/2/3/4/4/2 | 4/4/2/2/1/5 | 3/4/2/2/1/4 | 2/1/1/1/1/2 |
| $N_{nlp}$ | 2/2/2/2/2/2 | 3/3/7/7/7/2 | 32/28/26 /18/20/45 | 75/91/47 /42/34/122 | 67/83/69 /42/36/212 |
| $t_{post}$ (s) | - | - | - | 3600/3600/3600 /3600/1215/3600 | 0/0/0/0/0/0 |
| $N_{inf,post}$ | 0/0/0/0/0/0 | 2/2/3/4/4/2 | 4/4/2/2/1/5 | 3/4/2/2/1/4 | 0/0/0/0/0/0 |

Inf: infeasible solution; $N_{node}$: total number of nodes; $N_{inf}$: number of nodes with infeasibility; $N_{nlp}$: total number of NLP problems.
$t_{post}$: computational time for post check and refinement; $N_{inf,post}$: number of nodes with infeasibility after post check and refinement.

From Table 4 HCBB-RB and HCBB-FP generate very similar optimal solution with TAC of 4.904

M\$ year$^{-1}$ from the six initial points. While HCBB-RB requires 8%-21% less computational effort than HCBB-FP from the first, second, fifth and sixth initial points, it consumes 3%-46% more computational time from the third and fourth initial points due to the difficulties in solving NLP subproblem with stream flow rates appearing in the denominator to optimality when the flow rates approach zero[22,84]. From Table 4, HCBB-FP spends extremely long time (usually more than 1 hours) in post check and refinement procedure, while HCBB-RB takes negligible time because no refinement procedure is required after post check for those nodes flagged as infeasible.

Similar to Examples 1-3, some bypass efficiencies in the optimal solution from HCBB-RB are still fractional due to limited convergence tolerance or local optimum. We round them up and solve the resulting NLP problem again to generate the optimal profits of 4.901-4.905 M\$ year$^{-1}$, a relative change of around 0.02%. The final optimal design is shown in Fig. S10.

We also test GAMS/BRAON as the solver of each NLP subproblem in the customized B&B algorithm. It seems that the performance of the algorithms with GAMS/BARON as the NLP solver is worse than that of GAMS/CONOPT as the NLP solver. We will investigate the performance of the HCBB algorithms with GAMS/BARON as the NLP solver in the future work.

**Conclusion**

In this work, we proposed two homotopy continuation enhanced branch and bound (HCBB) algorithms namely HCBB-FP and HCBB-RB through using the homotopy continuation (HC) method to solve large-scale strongly nonconvex MINLP problems. During branch and bound (B&B) in both HCBB-FP and HCBB-RB, each NLP subproblem at a node was solved using the solution from its parent node as an initial point. While HCBB-FP solved a series of feasibility problems to gradually reach the feasible solution at each node, HCBB-RB solved a few optimality problems with gradually tightened bounds of the homotopy variable and reached the optimal solution. As a result, the HCBB-RB could terminate once the current optimum was larger than the upper bound of the MINLP problem during HC, which significantly reduced the computational effort required. A variable step length was adapted to effectively balance feasibility and computational efficiency. Several matching strategies were proposed to identify suitable step lengths based on historic information of the previous nodes, which reduced the effort in finding appropriate step lengths at latter similar nodes. To further improve solution quality, a post check and refinement procedure was proposed after B&B to revisit the nodes flagged as infeasible due to the minimum step length or the maximum number of iterations reached.

The computational results of solving process synthesis problems using rigorous unit operation models demonstrate that both HCBB-FP and HCBB-RB were able to generate the same best local optimum from different initial points, while all the other existing algorithms/solvers, including GAMS/SBB, GAMS/DICOPT, SRMINLP and our custom B&B algorithm failed to solve some examples or obtained much worse local optima. This is because the HC method can find optimal solutions at some nodes which were flagged as infeasible in GAMS/SBB, GAMS/DICOPT and our custom B&B algorithm, reducing the risk of missing a better solution. Therefore, HCBB-FP and HCBB-RB are more robust than the other existing algorithms. It is also demonstrated that HCBB-RB was superior to HCBB-FP in terms of the number of nodes flagged as infeasible, the computational time required and the solution quality. More importantly, HCBB-RB usually terminated the post check and refinement before entering the time-consuming refinement steps.

Although HCBB-FP and HCBB-RB have shown very good convergence performance for strongly nonlinear and nonconvex MINLP problems, they still have difficulties in solving NLP subproblems with stream flow rates appearing in the denominator to optimality when the flow rates approach zero[22,84], which inspires us to integrate the logic-based method in the future. Furthermore, it should be noted that the proposed algorithms cannot guarantee global optimality for nonconvex MINLP problems.

## Acknowledgements

The authors acknowledge financial support from China Scholarship Council – The University of Manchester Joint Scholarship (201809120005)

**Abbreviations**

| | |
|---|---|
| HCBB | Homotopy continuation enhanced branch and bound |
| HCBB-FP | HCBB algorithm with solving NLPFP |
| HCBB-RB | HCBB algorithm with solving NLPRB |
| BB | Branch and bound |
| HC | Homotopy continuation |
| OCFEM | Orthogonal collocation finite element method |
| HDA | hydrodealkylation |
| TAC | Total annualised cost |

**Conflict of Interest**

The authors claim no conflict of interest.

**Supporting Information**

Appendix A) The complete homotopy continuation enhanced branch and bound algorithm (HCBB); Appendix B) Nomenclature; Algorithm S1) Nonlinear B&B algorithm; Algorithm S2) Initialization of the whole process; Algorithm S3) Initialization of the reactor network using OCEFM models; Table S1) Rate constant, operating conditions and specifications in the reactor and distillation column for Example 1; Table S2) Stream temperature and vaporization enthalpy in preheater, column1 and column 2 for Example 1; Table S3) Component physical properties and Antoine coefficients for Example 1; Table S4) Comparative results from BB, HCBB-FP, and HCBB-RB for Example 1; Table S5) Rate constants, operating conditions, and specifications in the reactor and distillation column for Example 2; Table S6) Component physical properties and Antoine coefficients for Example 2; Table S7) Parameters in the Wilson equations for Example 2; Table S8) Comparative results from BB, HCBB-FP, and HCBB-RB for Example 2; Table S9) Feedstock and product specification and prices for Example 3; Table S10) Utility prices for Example 3; Table S11) Investment costs for Example 3; Table S12) Parameters for Antoine equation in Exmaple 3; Table S13) Parameters for ideal vapour heat capacity polynomial in Example 3; Table S14) Parameters for DIPPR liquid heat capacity polynomial in Example 3; Table S15) Comparative results from BB, HCBB-FP and HCBB-RB for Example 3; Table S16) Stoichiometric coefficients for HDA reactions; Table S17) Comparative results from BB, HCBB-FP and HCBB-RB for Example 4; Figure S1) Reaction pathways for benzene chlorination; Figure S2) Superstructure for the synthesis of benzene chlorination process; Figure S3) Optimal design of the benzene chlorination process; Figure S4) Reaction pathway for cyclohexane oxidation; Figure S5) Superstructure for the cyclohexane oxidation process; Figure S6) Optimal design of the cyclohexane oxidation process; Figure S7) Reaction pathways for the HDA; Figure S8) Superstructure of the HDA process; Figure S9) Superstructure for the synthesis of HDA process using differential reactor model; Figure S10) Optimal design of HDA process with differential reactor model.

Yingjie Ma and Jie Li[†]

Centre for Process Integration, Department of Chemical Engineering and Analytical Science,

School of Engineering, The University of Manchester, Manchester M13 9PL, UK

**List of Algorithms**

**Appendix A** The complete homotopy continuation enhanced branch and bound algorithm (HCBB); Algorithm S1) Standard branch and bound (BB) algorithm; Algorithm S2) Initialization of the whole process; Algorithm S3) Initialization of the reactor network using OCEFM models.

**List of Tables**

Table S1) Rate constant, operating conditions and specifications in the reactor and distillation column for Example 1; Table S2) Stream temperature and vaporization enthalpy in preheater, column1 and column 2 for Example 1; Table S3) Component physical properties and Antoine coefficients for Example 1; Table S4) Comparative results from BB, HCBB-FP, and HCBB-RB for Example 1; Table S5) Rate constants, operating conditions, and specifications in the reactor and distillation column for Example 2; Table S6) Component physical properties and Antoine coefficients for Example 2; Table S7) Parameters in the Wilson equations for Example 2; Table S8) Comparative results from BB, HCBB-FP, and HCBB-RB for Example 2; Table S9) Feedstock and product specification and prices for Example 3; Table S10) Utility prices for Example 3; Table S11) Investment costs for Example 3; Table S12) Parameters for Antoine equation in Exmaple 3; Table S13) Parameters for ideal vapour heat capacity polynomial in Example 3; Table S14) Parameters for DIPPR liquid heat capacity polynomial in Example 3; Table S15) Comparative results from BB, HCBB-FP and HCBB-RB for Example 3; Table S16) Stoichiometric coefficients for hydrodealkylation (H) reactions; Table S17) Comparative results from BB, HCBB-FP and HCBB-RB for Example 4.

---

[†] Corresponding author: Jie Li (Email address: jie.li-2@manchester.ac.uk). Tel: + 44 (0) 161 306 8622. Fax: +44 (0) 161 236 7439

**List of Figures**

Figure S1) Reaction pathways for benzene chlorination; Figure S2) Superstructure for the synthesis of benzene chlorination process; Figure S3) Optimal design of the benzene chlorination process; Figure S4) Reaction pathway for cyclohexane oxidation; Figure S5) Superstructure for the cyclohexane oxidation process; Figure S6) Optimal design of the cyclohexane oxidation process; Figure S7) Reaction pathways for the HDA; Figure S8) Superstructure of the HDA process; Figure S9) Superstructure for the synthesis of HDA process using differential reactor model; Figure S10) Optimal design of HDA process with differential reactor model where the unit of flow rate is kmol h$^{-1}$ and A, B, C, D and E represent $H_2$, $CH_4$, $C_6H_6$, $C_7H_8$ and $C_{12}H_{10}$ respectively.

## Appendix A The complete homotopy continuation enhanced branch and bound algorithm

---

**Algorithm A1**: HCBB algorithm with post check and refinement

---

    **Data**: $f$, $\mathbf{g}$, $\mathbf{h}$, $\mathbf{x}^0$, $\mathbf{y}^0$, $I$

1  **initialization**: $i \leftarrow 0$, $\mathcal{S}_F^i \leftarrow \varnothing$, $\mathcal{S}_R^i \leftarrow \mathcal{S}$, $\mathcal{L} \leftarrow \{n_0\}$, $f^{ub} \leftarrow +\infty$, incumbent solution $\leftarrow$ None, $\mathcal{I} \leftarrow \varnothing$;

2  **while** $\mathcal{L} \neq \emptyset$ **do**

3      select a node $n_i$ from $\mathcal{L}$;

4      solve the current subproblem with **Algorithm 2**;

5      **if** the solution is infeasible **then**

6          put the current node $n_i$ to $\mathcal{I}$, record $\mathbf{y}_F^{i,last}$, the last feasible or optimal solution $(\mathbf{x}^{i,last}, \mathbf{y}_R^{i,last}, f^{i,last})$ and the corresponding $t^{i,last}$ and $\Delta t^{i,last}$ during HC;

7      prune the branch, update the incumbent or branch the node as with lines 5-16 of the **Algorithm S1**;

8  **while** $\mathcal{I} \neq \emptyset$ **do**

9      select a node $i$ from $\mathcal{I}$ and let $\mathcal{L} \leftarrow \{i\}$, set $f^{ub} \leftarrow f^*$;

10    apply post check and refinement procedure (its lines 3-19) to the node

11    **if** the node is fathomed **then**

12       prune the branch;

13    **else**

14       go to the line 2 of the current algorithm and conduct B&B again;

---

**Appendix B Nomenclature**

**Sets**

$\mathcal{S}$: Set of indices of binary variables;

$\mathcal{I}$: Set of infeasible nodes in HCBB algorithms;

$\mathcal{J}^i$: Set of the nodes $n_j$ ($j < i$) having the same branch variable as with the current node $n_i$;

$\mathcal{L}$: Set of all the nodes during branch and bound;

$\mathcal{S}_F^i$: Set of indices of fixed binary variables at node $n_i$;

$\mathcal{S}_R^i$: Set of indices of relaxed binary variables at node $n_i$;

$\mathcal{T}$: Collection of all the $\Delta t^{v,j}$ and $t^{v,j}$;

$\mathbb{R}^n$: $n$-dimension real number set;

$\mathbb{B}_{0,1}^{i,R}$: $m_R^i$-dimension rectangular region between 0 and 1 at node $n_i$;

$\mathbb{B}_{0,1}^{i,F}$: $m_F^i$-dimension rectangular region between 0 and 1 at node $n_i$;

$\mathbb{Z}^{m_F^i}$: $m_F^i$-dimension Set of binary variables at node $n_i$;

**Indices**

$s$: Index of a binary variable;

$i, j$: Index of a node;

$bpi$: Index of the branch variable at node $n_{pi}$;

$pi$: Index of the parent node of node $n_i$;

$d$: Index of a fixed binary variable;

$l$: Index of a relaxed binary variable;

**Parameters**

$m$: Dimension of binary variables;

$m_F^i$: Dimension of fixed binary variables at node $n_i$;

$m_R^i$: Dimension of relaxed binary variables at node $n_i$;

$n$: Dimension of continuous variables;

$n_i, n_j$: The nodes to be explored at the $i^{th}$, $j^{th}$ iterations respectively during branch and bound;

$n_{pi}, n_{pj}$: The parent nodes of nodes $n_i, n_j$ respectively during branch and bound;

$N$: Maximum allowable number of steps for homotopy continuation (HC);

$N^{post}$: Maximum allowable number of steps for HC during post check and refinement;

$t$: Homotopy parameter;

$t^{i,last}$: Homotopy parameter corresponding to the last feasible or optimal solution during HC at an infeasible node $n_i$.

$\Delta t$: Step length of the homotopy parameter during homotopy continuation;

$\Delta t^{i,last}$: The last step length of the homotopy parameter during homotopy continuation at an infeasible node $n_i$;

$\Delta t^{v,j}$: the $v^{th}$ step length of homotopy parameter leading to a feasible or optimal solution during homotopy continuation of the node $n_j$;

$\Delta t_{min}$: Minimum step length of the homotopy parameter $t$;

$\Delta t_{min}^{post}$: Minimum step length of the homotopy parameter $t$ during post check and refinement;

$\Delta y_b^{i,j}$: Difference between the optimal values of the branching variables at nodes $n_{pi}$ and $n_{pj}$;

$\delta$: Threshold to indicate whether $\Delta y_b^{i,j}$ is small enough;

$v$: Iteration number during the homotopy continuation process.

$t^{v,j}$: The $v^{th}$ homotopy parameter leading to a feasible or optimal solution during homotopy continuation of the node $n_j$;

**Variables**

$f^*$: An optimum of the MINLP problem generated from the HCBB algorithm;

$f^{i,*}$: An optimum of the NLP subproblem at node $n_i$;

$f^{v,*}$: An optimum of the NLP subproblem at a homotopy continuation iteration $v$;

$f^{i,last}$: The last feasible or optimal solution in homotopy continuation at an infeasible node $n_i$;

$f^{ub}$: The upper bound of the MINLP problem;

$\mathbf{x}$: Continuous variables to be optimized;

$\mathbf{x}^{i,*}$, $\mathbf{x}^{pi,*}$: Optimal values of $\mathbf{x}$ from an NLP subproblem at node $n_i$ and $n_{pi}$ respectively;

$\mathbf{x}^{i,last}$: The values of continuous variables in the last feasible or optimal solution during homotopy continuation at an infeasible node $n_i$;

$\mathbf{y}$: Binary variables to be optimized;

$y_{bpi}^i$, $y_{bpi}^j$: Required value of the binary variable with index $bpi$ at node $n_i$ and $n_j$ respectively;

$\tilde{y}_{bpi}^i$: Homptopy variable during HC;

$y_{bpi}^{pi,*}$, $y_{bpj}^{pj,*}$: Optimal values of the binary variable for branching at nodes $n_{pi}$ and $n_{pj}$ respectively;

$\mathbf{y}_F^i$, $\mathbf{y}_F^{pi}$: Binary variables whose values are known and fixed at node $n_i$ and node $n_{pi}$ respectively;

$\mathbf{y}_{F,s}^i\big|_{s \neq bpi}$: A vector of the binary variables in $\mathbf{y}_F^i$ at node $n_i$ whose indices are not $bpi$;

$\tilde{\mathbf{y}}_F^i$: The intermediate values of fixed binary variables $y_F^i$;

$\mathbf{y}_R^i$: relaxed binary variables to be optimized at node $n_i$;

$\mathbf{y}_R^{i,*}$, $\mathbf{y}_R^{pi,*}$: Optimal values of the relaxed binary variables of an NLP subproblem at node $n_i$ and $n_{pi}$ respectively;

$\mathbf{y}_R^{i,last}$: the values of relaxed binary variables in the last feasible or optimal solution during HC at an infeasible node $n_i$;

# 1 Nonlinear branch and bound algorithm

A nonlinear B&B algorithm[1] for a general MINLP problem is presented in **Algorithm S1**.

---

**Algorithm S1**: Nonlinear B&B algorithm

---

**Data:** $f$, $g$, $h$, $I$

1  **initialization**: $i \leftarrow 0$, $\mathcal{S}_F^i \leftarrow \varnothing$, $\mathcal{S}_R^i \leftarrow \mathcal{S}$, $\mathcal{L} \leftarrow \{n_0\}$, $f^{ub} \leftarrow +\infty$, incumbent solution $\leftarrow$ None;

2  **while** $i < I$ and $\mathcal{L} \neq \emptyset$ **do**

3       select a node $n_i$ from $\mathcal{L}$;

4       fix $\mathbf{y}_F^i$ at specified binary values, relax $\mathbf{y}_R^i$ to be continuous values between 0 and 1, and then solve the resulting NLP subproblem. The optimal solution has $\mathbf{x}^{i,*}$, $\mathbf{y}_R^{i,*}$ and $f^{i,*}$;

5       **if** the current NLP subproblem is infeasible **then**

6           prune the current branch;

7       **else**

8           **if** $f^{i,*} \geq f^{ub}$ **then**

9               prune the current branch ;

10           **else if** all components of $\mathbf{y}^{i,*}$ are binaries **then**

11               incumbent solution$\leftarrow (\mathbf{x}^{i,*}, \mathbf{y}_R^{i,*}, f^{i,*})$, $f^{ub} \leftarrow f^i$;

12               prune the current branch;

13           **else**

14               select a variable with a fractional value in $\mathbf{y}_R^{i,*}$; generate two nodes through fixing the variable at 0 and 1 respectively; put the two nodes in $\mathcal{L}$;

15       remove $n_i$ from $\mathcal{L}$; $i \leftarrow i + 1$;

16  **if** the incumbent solution is not None **then**

17       The optimal integer solution is found;

18  **else**

19       No feasible integer solution found for the problem;

---

# 2 Initialization methods

The initialization method from Ma et al.[2] is presented in Algorithm S2.

| | **Algorithm S2**: Initialization of the whole process |
|---|---|
| 1 | solve the reactor network model without economic evaluation equations to get a feasible solution using BARON; |
| 2 | solve economic evaluation equations using CONOPT; |
| 3 | solve the reactor network model with the above solution as an initial point and get the optimal solution using CONOPT; |
| 4 | simulate or optimize the distillation system with the effluent from the reactor network as a feed to the distillation system; |
| 5 | optimize the reaction-separation-recycle process for a feasible solution with the above solutions from reactor network and distillation system as an initial point; |

In step 1, the global optimizer BARON is used to find a feasible solution for the reactor network described by algebraic equations within acceptable time as the scale of such model is small. However, it cannot generate a feasible solution for the reactor network described by differential equations (see Example 4) due to a large number of equations in such model after discretization using the orthogonal collocation method[3]. The initialization method from Ma et al.[2] is extended to tackle such a more challenging case, as shown in the Algorithm S3 where $M$ is the number of collocation points required, $k$ is the number of collocation points in the current problem, and $\bar{k}$ is the number of collocation points in the latest feasible solution.

In the Algorithm S3, Line 8 is to increase the number of collocation points when the solution is feasible, while line 10 is to decrease the number of collocation points when it is infeasible. In line 11, the results from the converged solution will be interpolated to provide good initial values for the model using more collocation points. In lines 5-12, when we increase the number of collocation points gradually, only a feasible solution is required at each iteration, while finally, the optimal solution will be searched in line 13. In Algorithm S3, the results from GAMS communicate information with Python through the gdx-pandas package[5] and Pandas[6] is used to transform data.

**Algorithm S3**: Initialization of the reactor network using OCFEM models

    **Data**: $M$;

1   **initialization**: $\bar{k} \leftarrow 0$, the number of collocation points for the reactor model $k \leftarrow 1$;

2   solve the reactor network model without economic evaluation equations to get a feasible solution using BARON;

3   solve economic evaluation equations using CONOPT;

4   specify the solution from previous two parts as the initial point;

5   **while** $k < M$ and $k \neq \bar{k}$ **do**

6      solve the reactor network model with $k$ collocation points using CONOPT to get a feasible solution from a given initial point;

7      **if** feasible solution is got **then**

8          $\bar{k} \leftarrow k$, $k \leftarrow \min(2k, M)$;

9      **else**

10          $k \leftarrow round(\frac{k+\bar{k}}{2})$;

11      interpolate the variable values for the model with $k$ collocation points using the Lagrange polynomial[4];

12      specify the interpolated results as initial values;

13   solve the reactor network model with the feasible solution as an initial point to get the optimal solution using CONOPT

## 3 Example 1

### 3.1 Reactions and parameters

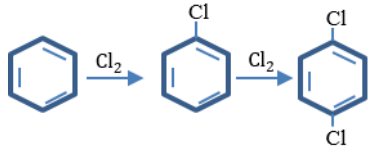

Figure S1 Reaction pathway for benzene chlorination

All reactions are first order. The reaction rate is equal to the kinetic constants multiplied by the molar concentration of benzene and chlorobenzene respectively for the two reactions. The rate constants, operating conditions and specifications of the reactor and distillation column are given in Table S1.

Table S1 Rate constant, operating conditions and specifications in the reactor and distillation column for example 1

| Reactor | |
|---|---|
| $k_1$ | 0.412 h$^{-1}$ |
| $k_2$ | 0.055 h$^{-1}$ |
| Reaction temperature | 30 $^{\circ}$C |
| Reaction pressure | 1 atm |
| Chlorobenze product rate | 50 kmol h$^{-1}$ |
| **Distillation column** | |
| Operating pressure | 1 atm |
| Upper bound of the number of trays | 25 |
| Initial feed location (top to bottom) | 13 |

The enthalpy is calculated by Eqs. (S1-S2).

$$h_L^i = h_{f,L}^i + Cp_L^i(T - T_0), \tag{S1}$$

$$h_V^i = h_{f,G}^i + Cp_V^i(T - T_0), \tag{S2}$$

where $h_{f,L}^i$ and $h_{f,V}^i$ are the liquid and vapour standard formation enthalpies respectively. $Cp_L^i$ and $Cp_V^i$ are mass specific heat capacities for liquid and vapour respectively. $T_0$ is the reference temperature, which is 273.15 K. $T$ (K) is the process stream temperature.

Table S2 Stream temperature and vaporization enthalpy in preheater, column1 and column 2 for example 1

|  | Preheater | Column 1 | Column 2 |
|---|---|---|---|
| Temperature (K) | 425 | 450 | 474.5 |
| Vaporization enthalpy (kJ kg⁻¹) | 2107.42 | 2021.4 | 1933.19 |

Table S3 Component physical properties and Antoine coefficients for Example 1

| Item | Benzene | Chlorobenzene | Dichlorobenzene |
|---|---|---|---|
| Molecular weight (kg kmol⁻¹) | 78.11 | 112.56 | 147.01 |
| Density (kg m⁻³) | 876 | 1110 | 1300 |
| Vapour standard formation enthalpy (kJ kmol⁻¹) | 82900 | 520000 | 33000 |
| Liquid standard formation enthalpy (kJ kmol⁻¹) | 48700 | 11100 | −17400 |
| Vapour heat capacity (kJ kmol⁻¹ K⁻¹) | 94 | 108 | 118 |
| Liquid heat capacity (kJ kmol⁻¹ K⁻¹) | 134 | 161 | 193 |
| Antoine coefficient A | 13.78 | 13.86 | 14.26 |
| Antoine coefficient B | 2726.81 | 3174.78 | 3798.20 |
| Antoine coefficient C | −55.58 | −61.45 | −59.83 |

The saturate pressure of each component is computed using the Antoine equation, as Eq. (S3)

$$P_{sat}^i = \exp(A^i - \frac{B^i}{C^i + T - T_0}), \tag{S3}$$

where $A^i$, $B^i$ and $C^i$ are Antoine coefficients for each component, which are shown in Table S2.

### 3.2 Superstructure for the synthesis of benzene chlorination process

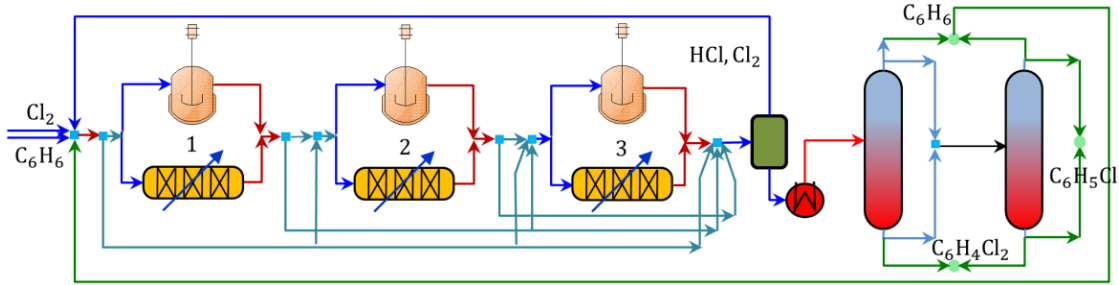

Figure S2 Superstructure for the synthesis of benzene chlorination process

### 3.3 Optimisation results

Table S4 Comparative results from BB, HCBB-FP, and HCBB-RB for Example 1

| Item | BB | HCBB-FP | HCBB-RB |
|---|---|---|---|
| Initial point 1 | | | |
| TAC | 0.613 | 0.613 | 0.613 |
| TAC after rounding | 0.613 | 0.613 | 0.613 |
| Relative change | 0.00% | 0.00% | 0.00% |
| Initial point 2 | | | |
| TAC | 0.614 | 0.614 | 0.614 |
| TAC after rounding | 0.614 | 0.614 | 0.614 |
| Relative change | 0.00% | 0.00% | 0.00% |
| Initial point 3 | | | |
| TAC | 0.614 | 0.614 | 0.614 |
| TAC after rounding | 0.614 | 0.614 | 0.614 |
| Relative change | 0.00% | 0.00% | 0.00% |

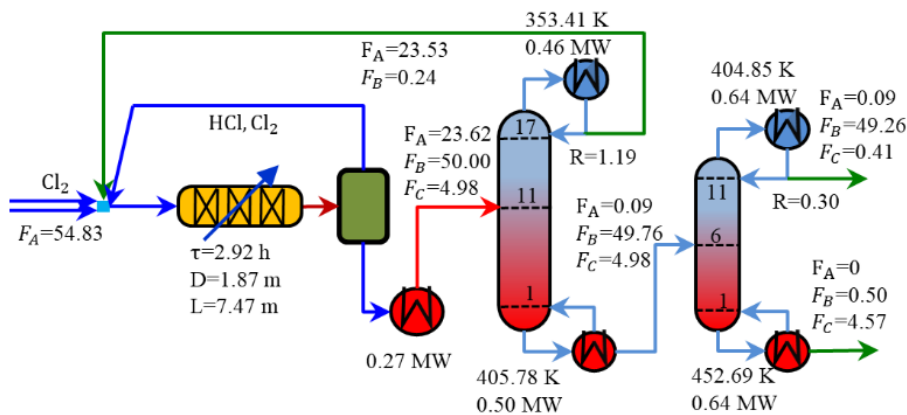

Figure S3 Optimal design of the benzene chlorination process

## 4   Example 2

### 4.1   Reactions and parameters

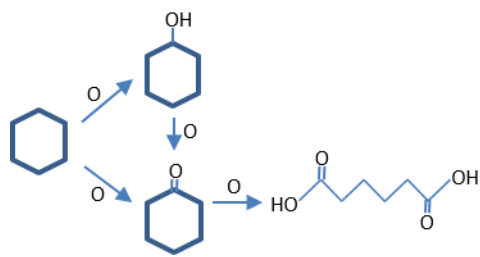

Figure S4 Reaction pathway for cyclohexane oxidation.

All reactions are first order. The reaction rates are equal to the kinetic constants multiplied by the molar concentration of cyclohexane, cyclohexane, cyclohexanol and cyclohexanone respectively for the four reactions. Some parameters and specifications in the distillation columns are given in Table S5.

Table S5 Rate constants, operating conditions, and specifications in the reactor and distillation column for Example 2

| **Reactor** | |
| --- | --- |
| $k_1$ | 0.12 h$^{-1}$ |
| $k_2$ | 0.032 h$^{-1}$ |
| $k_3$ | 2.85 h$^{-1}$ |
| $k_4$ | 2.30 h$^{-1}$ |
| Reaction temperature | 160 ºC |
| Reaction pressure | 10 atm |
| Chlorobenzene product rate | 50 kmol h$^{-1}$ |
| **Distillation column** | |
| Operating pressure | 1 atm |
| Upper bound of the number of trays | 25 |
| Initial feed location (top to bottom) | 13 |

Table S6 Component physical properties and Antoine coefficients for Example 2

| | Cyclohexane | Cyclohexanone | Cyclohexanol | Adipic acid |
| --- | --- | --- | --- | --- |
| Molecular weight (kg kmol$^{-1}$) | 84.2 | 98.1 | 100.2 | 337.5 |
| Density (kg m$^{-3}$) | 778 | 948 | 962 | 1360 |
| Vapour standard formation enthalpy (kJ kmol$^{-1}$) | –123100 | –231100 | –29000 | –659400 |
| Liquid standard formation enthalpy (kJ kmol$^{-1}$) | –156400 | –276100 | –352000 | –741300 |
| Vapour heat capacity (kJ kmol$^{-1}$ K$^{-1}$) | 105 | 161 | 176 | 245 |
| Liquid heat capacity (kJ kmol$^{-1}$ K$^{-1}$) | 156 | 177 | 214 | 330 |
| Antoine coefficient A | 13.74 | 12.39 | 14.05 | 17.21 |
| Antoine coefficient B | 2766.63 | 2101.96 | 3443.70 | 6477.32 |
| Antoine coefficient C | -50.50 | -164.02 | -63.59 | -95.95 |

Similar to Example 1, the enthalpy and saturate pressure of each component are calculated by Eqs. (S1-S3). However, Wilson equation is used to calculate liquid activity coefficients instead of assuming ideal liquid mixture, which is shown in Eq. S4

$$ln\gamma_i = 1 - \ln(\sum_{k=1}^{4} x_k \Lambda_{ik}) - \sum_{l=1}^{4} \frac{x_l \Lambda_{li}}{\sum_{k=1}^{4} x_k \Lambda_{kl}}, \tag{S4}$$

where $\gamma_i$ is the component activity coefficient, $x_k$ and $x_l$ are molar compositions, and $\Lambda_{ik}$, $\Lambda_{li}$ and $\Lambda_{kl}$ are calculated from Eq. S5.

$$\Lambda_{ij} = \frac{MV_j}{MV_i} \exp(\frac{\lambda_{ij} - \lambda_{ii}}{RT}). \tag{S5}$$

The binary coefficients $\lambda_{ij}$ and the molar volume $MV_i$ are provided in Table S7.

Table S7 Parameters in the Wilson equations and molar volume for Example 2

| $i$ \ $j$ | | Cyclohexane | Cyclohexanone | Cyclohexanol | Adipic acid |
|---|---|---|---|---|---|
| $\lambda_{ij}$ | Cyclohexane | 0 | 310.99 | -409.97 | 2576.96 |
| | Cyclohexanone | 6074.00 | 0 | 1445.72 | -848.27 |
| | Cyclohexanol | 2560.35 | 357.55 | 0 | -1409.18 |
| | Adipic acid | 9809.50 | 523.97 | 4965.28 | 0 |
| MV (m³ kmol⁻¹) | | 0.10816 | 0.10407 | 0.10356 | 0.10746 |

The steam temperature is defined to be 30 K higher than the bubble point of the mixture in the reboiler (LMTD = 30), that is $(T_{reb} + 30 - 273.15)$ ℃, where $T_{reb}$ is the temperature in the reboiler, which can be got from MESH equations. Furthermore, the steam latent heat $\Delta H_{vap}$ is calculated as follows:

$$\Delta H_{vap} = -0.005 \cdot (T_{reb} - 243.15)^2 - 1.710 \cdot (T_{reb} - 243.15) + 2483.050, \tag{S6}$$

### 4.2 Superstructure for the synthesis of cyclohexane oxidation process

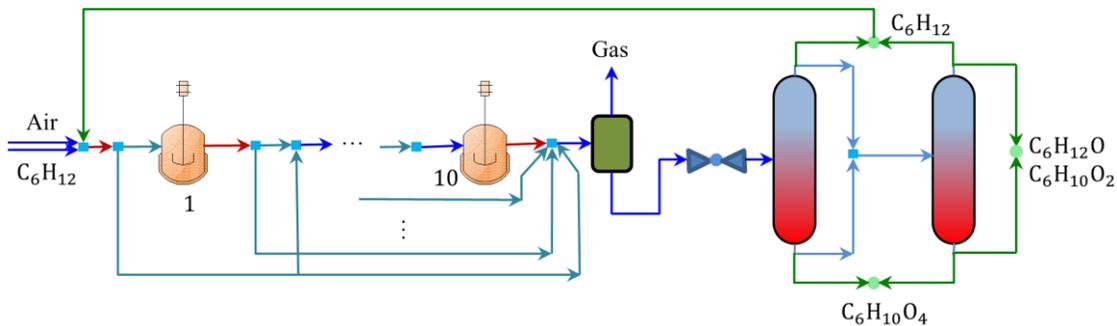

Figure S5 Superstructure for the cyclohexane oxidation process

### 4.3 Computational results

Table S8 Comparative results from BB, HCBB-FP, and HCBB-RB for Example 2

|  | BB | HCBB-FP | HCBB-RB |
|---|---|---|---|
| TAC | Inf | 2.146 | 2.118 |
| TAC after rounding | - | 2.148 | 2.129 |
| Relative change | - | 0.1% | 0.5% |
| TAC | Inf | 2.146 | 2.118 |
| TAC after rounding | - | 2.148 | 2.120 |
| Relative change | - | 0.1% | 0.1% |
| TAC | 3.572 | 2.146 | 2.118 |
| TAC after rounding | 3.591 | 2.152 | 2.129 |
| Relative change | 0.5% | 0.35 | 0.5% |

Inf: infeasible.

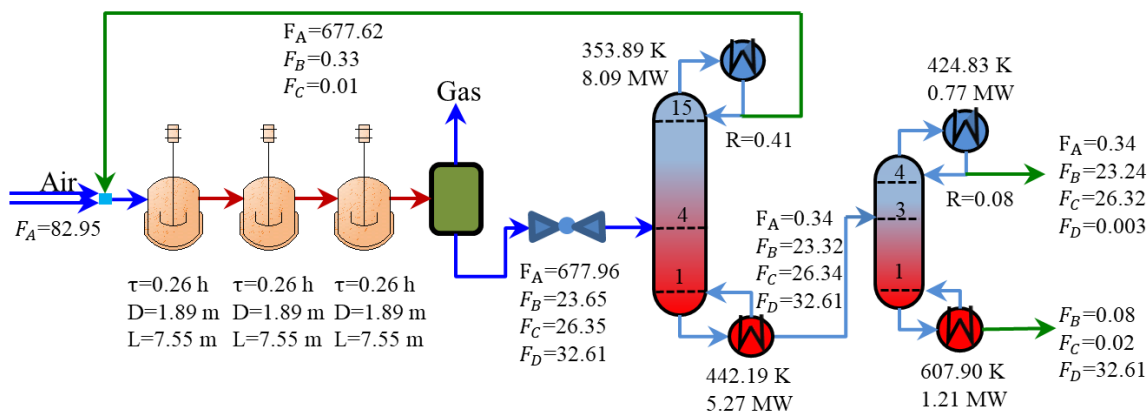

Figure S6 Optimal design of the cyclohexane oxidation process

## 5 Example 3

### 5.1 Reactions

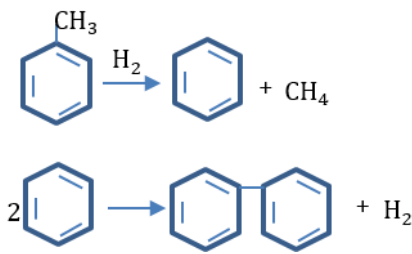

Figure S7 Reaction pathways for the HDA

The kinetic equations for an isothermal PFR are given as follows,

$$1 - X = \frac{1}{\left(1 + 0.372kVF_{in}^{tol\,0.5}F_{in}^{tot\,-1.5}\right)^2}, \tag{S7}$$

$$1 - Y = 0.0036(1 - X)^{-1.544}, \tag{S8}$$

$$k = 6.3 \times 10^{10} \exp\frac{-26167}{RT}, \tag{S9}$$

where $X$ is toluene conversion, $k$ (kmol m³ s⁻¹) is the reaction rate constant, $V$ (m³) is reactor volume, $F_{in}^{tol}$ (kmol h⁻¹) is the molar flow rate of toluene in the feed stream, $Y$ is the selectivity of benzene, $R$ (kJ kmol K⁻¹) is ideal gas constant and $T$ (K) is reaction temperature. An adiabatic reactor is approximated by an isothermal reactor with the reaction temperature being the average value of the inlet and outlet temperatures.

### 5.2 Superstructure for the synthesis of HDA process

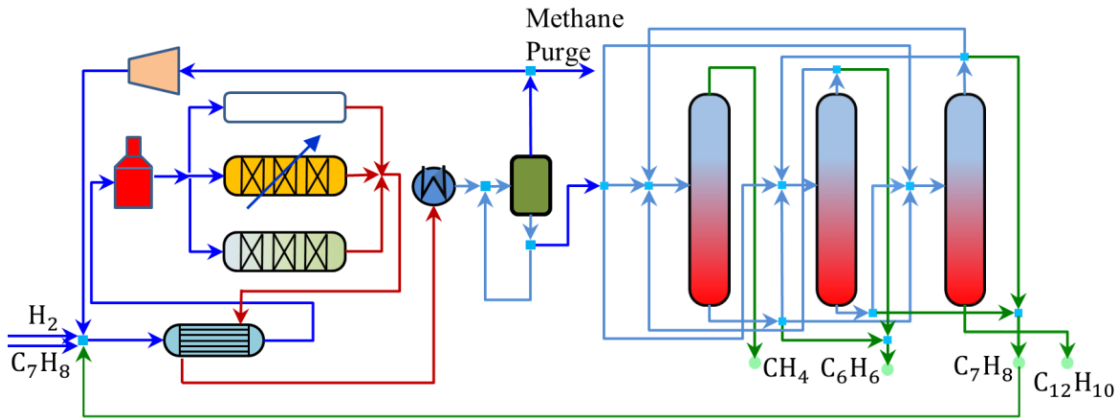

Figure S8 Superstructure of the HDA process

### 5.3 Cost model for HDA

Table S9 Feedstock and product specification and prices for Example 3

|  | Specification | Price ($ kmol⁻¹) |
|---|---|---|
| Hydrogen feed | 95% hydrogen, 5% methane | 2.50 |
| Toluene feed | 100% toluene | 14.00 |
| Benzene product | ≥ 99.97% benzene | 19.90 |
| Diphenyl product | - | 11.84 |
| Hydrogen purge | Heating value | 1.08 |
| Methane purge | Heating value | 3.37 |

<div align="center">Table S10 Utility prices for Example 3</div>

|  | Price |
|---|---|
| Electricity | 0.04 \$ (kW·h)$^{-1}$ |
| Steam | 8.0 \$ GJ$^{-1}$ |
| Water | 0.7 \$ GJ$^{-1}$ |
| Fuel | 4.0 \$ GJ$^{-1}$ |

<div align="center">Table S11 Investment costs for Example 3</div>

|  | Fixed-charge cost (k\$ year$^{-1}$) | Variable-charge cost (k\$ year$^{-1}$) |
|---|---|---|
| Compressor | 7.155 | 0.815 × brake horsepower (kW) |
| Stabilizing column | 1.126 | 0.375 × number of stages |
| Benzene column | 16.3 | 1.55 × number of stages |
| Toluene column | 3.90 | 1.12 × number of stages |
| Furnace | 6.20 | 1.172 × heat duty (103 GJ year$^{-1}$) |
| Adiabatic reactor | 74.3 | 1.257 × reactor volume (m$^3$) |
| Isothermal reactor | 92.875 | 1.571 × reactor volume (m$^3$) |

The compressor brake horsepower $Power$ is calcuated from Eq. (S10), assuming isentropic compression.

$$Power = \frac{\eta}{eff}\theta\frac{V_{in}}{3600}T_{in}(\alpha - 1), \tag{S10}$$

$$\theta = \frac{CpCv}{CpCv-1}, \tag{S11}$$

$$\alpha^{\theta} = \frac{P_{out}}{P_{in}} = \frac{T_{out}}{T_{in}}, \tag{S12}$$

where $\eta$ and $eff$ are two constants, which are 0.3665 and 0.75 respectively accroding to the GAMS file provided by Kocis and Grossmann.[7] $\alpha$ is the compression ratio. $V_{in}$ (kmol h$^{-1}$) and $T_{in}$ (K) are the flow rate and temperature of the inlet stream of the compressor. $P_{in}$ (Pa) and $P_{out}$ (Pa) are the inlet and outlet pressures of the compressor respectively and $T_{out}$ (K) is the outlet temperature. $CpCv$ is the ratio of the heat capacity at constant pressure and heat capacity at constant volume, which is 1.3 in this work.

The heat duty of the furnace can be got from the heat balance between the inlet and outlet stream easily, while the models of distillation columns and reactors are shown in the main text.

### 5.4 Saturate vapor pressure calculation

The Antoine equation is used to calculate the saturate vapor pressure, as shown in Eq. (S13).

$$\log(750 \cdot P_i^s) = A_i - \frac{B_i}{C_i + T} \qquad\qquad \forall i \qquad\qquad (S13)$$

where the unit of saturate pressure $P_i^s$ is bar and the unit of temperature $T$ is K. The parameters are from Kocis and Grossmann[8] and shown in the following Table S12.

Table S12 Parameters in the Antoine equation in Example 3

|            | A       | B       | C      |
|------------|---------|---------|--------|
| $H_2$      | 13.6333 | 164.9   | 3.19   |
| $CH_4$     | 15.2243 | 897.84  | -7.16  |
| $C_6H_6$   | 15.9008 | 2788.51 | -52.36 |
| $C_7H_8$   | 16.0137 | 3096.52 | -53.67 |
| $C_{12}H_{10}$ | 16.6832 | 4602.23 | -70.42 |

### 5.5 Ideal gas heat capacity

The ideal gas heat capacity of component $i$, $C_p^{i,g}$ is calculated according to Aspen ideal gas heat capacity polynomial, as shown in Eq. (S14).

$$C_p^{i,g} = C_{g1} + C_{g2}T + C_{g3}T^2 + C_{g4}T^3 + C_{g5}T^4 + C_{g6}T^5 \qquad\qquad \forall i \qquad\qquad (S14)$$

Here, the unit of ideal gas heat capacity $C_p^{i,g}$ is J kmol$^{-1}$ K$^{-1}$ and the unit of temperature is K. The parameters from Aspen Properties V8.8[9] are shown in Table S13.

Table S13 Parameters for ideal vapor heat capacity polynomial in Example 3

|            | $C_{g1}$   | $C_{g2}$ | $C_{g3}$ | $C_{g4}$             | $C_{g5}$              | $C_{g6}$               |
|------------|------------|----------|----------|----------------------|-----------------------|------------------------|
| $H_2$      | 25555.306  | 14.821   | -0.023   | $1.768\times10^{-5}$ | $-5.935\times10^{-9}$ | $7.135\times10^{-13}$ |
| $CH_4$     | 19250.906  | 52.126   | 0.012    | $-1.132\times10^{-5}$ | 0                     | 0                      |
| $C_6H_6$   | -33917.267 | 474.364  | -0.302   | $7.130\times10^{-5}$ | 0                     | 0                      |
| $C_7H_8$   | -24354.616 | 512.464  | -0.277   | $4.911\times10^{-5}$ | 0                     | 0                      |
| $C_{12}H_{10}$ | -97066.771 | 1105.734 | -0.886   | 0.0003               | 0                     | 0                      |

### 5.6 Liquid heat capacity

The liquid heat capacity of each component $C_p^{i,l}$ is calculated according to DIPPR liquid heat capacity correlation, as shown in Eq. (S15).

$$C_p^{i,l} = C_{l1}^i + C_{l2}^i T + C_{l3}^i T^2 + C_{l4}^i T^3 + C_{l5}^i T^4 \qquad \forall i \qquad \text{(S15)}$$

where the unit of liquid heat capacity $C_p^{i,l}$ is J kmol$^{-1}$ K$^{-1}$ and the unit of temperature is K. The parameters from Aspen Properties V8.8 are shown in Table S14. Note that the liquid enthalpies of H$_2$ and CH$_4$ are calculated according to their vapor enthalpies, so the parameters for them are not shown in Table S14.

Table S14 Parameters used in DIPPR liquid heat capacity polynomial in Example 3

|  | $C_{l1}$ | $C_{l2}$ | $C_{l3}$ | $C_{l4}$ | $C_{l5}$ |
|---|---|---|---|---|---|
| C$_6$H$_6$ | 129440 | –169.5 | 0.648 | 0 | 0 |
| C$_7$H$_8$ | 140140 | –152.3 | 0.695 | 0 | 0 |
| C$_{12}$H$_{10}$ | 121770 | 429.3 | 0 | 0 | 0 |

### 5.7 Specific enthalpies of vapor and lquid

The vapor specific enthalpy is calculated from Eq. (S16),

$$h_i^g = h_i^{g,0} + \int_{T_0}^{T} C_{p,i}^{*,g} dT \qquad \forall i \qquad \text{(S16)}$$

where $h_i^{g,0}$ and $h_i^g$ are the vapor specific enthalpies of component $i$ at temperatures $T_0$ and $T$ respectively. $T_0$ is a referring temperature, which is 300 K. When the ideal gas heat capacity is calculated from Eq. (S15), Eq. (S16) can be integrated analytically, leading to Eq. (S17),

$$h_i^g = h_i^{g,0} + C_{g1}(T - T_0) + \tfrac{1}{2}C_{g2}(T^2 - T_0^2) + \tfrac{1}{3}C_{g3}(T^3 - T_0^3) + \tfrac{1}{4}C_{g4}(T^4 - T_0^4) +$$
$$\tfrac{1}{5}C_{g5}(T^5 - T_0^5) + \tfrac{1}{6}C_{g6}(T^6 - T_0^6) \qquad \forall i \qquad \text{(S17)}$$

The liquid specific enthalpy $h_i^l$ is calculated from Eq. (S18),

$$h_i^l = h_i^{l,0} + \int_{T_0}^{T} C_{p,i}^{*,l} dT \qquad \forall i \qquad \text{(S18)}$$

Here, $h_i^{l,0}$ and $h_i^l$ are the liquid specific enthalpies of component $i$ at temperatures $T_0$ and $T$ respectively. When applying DIPPR liquid heat capacity polynomial, Eq. (S18) can be integrated analytically, resulting in Eq. (S19),

$$h_i^l = h_i^{l,0} + C_{l1}(T - T_0) + \tfrac{1}{2}C_{l2}(T^2 - T_0^2) + \tfrac{1}{3}C_{l3}(T^3 - T_0^3) + \tfrac{1}{4}C_{l4}(T^4 - T_0^4) +$$
$$\tfrac{1}{5}C_{l5}(T^5 - T_0^5) \qquad \forall i \qquad \text{(S19)}$$

In Eqs. (S17) and (S19), both $h_i^{g,0}$ and $h_i^{l,0}$ are got from Aspen Plus V8.8 under reference temperature $T_0$, so they have considered the component formation enthalpies. Hence, $h_i^g$ and $h_i^l$ also include the component formation enthalpies.

### 5.8 Computational results

Table S15 Comparative results from BB, HCBB-FP and HCBB-RB for Example 3

| Result | BB | HCBB-FP | HCBB-RB |
|---|---|---|---|
| | | Init 1 | |
| Profit | 4.956 | 4.956 | 4.956 |
| Profit after rounding | 4.954 | 4.954 | 4.954 |
| Relative change | 0.04% | 0.04% | 0.04% |
| | | Init 2 | |
| Profit | 4.954 | 4.954 | 4.954 |
| Profit after rounding | 4.952 | 4.952 | 4.952 |
| Relative change | 0.04% | 0.04% | 0.04% |
| | | Init 3 | |
| Profit | 4.956 | 4.956 | 4.956 |
| Profit after rounding | 4.954 | 4.954 | 4.954 |
| Relative change | 0.04% | 0.04% | 0.04% |
| | | Init 4 | |
| Profit | 4.958 | 4.958 | 4.958 |
| Profit after rounding | 4.956 | 4.956 | 4.956 |
| Relative change | 0.04% | 0.04% | 0.04% |
| | | Init 5 | |
| Profit | 4.955 | 4.955 | 4.955 |
| Profit after rounding | 4.953 | 4.953 | 4.953 |
| Relative change | 0.04% | 0.04% | 0.04% |
| | | Init 6 | |
| Profit | 4.958 | 4.958 | 4.958 |
| Profit after rounding | 4.956 | 4.956 | 4.956 |
| Relative change | 0.04% | 0.04% | 0.04% |

## 6 Example 4

### 6.1 Kinetic equation for HDA

The kinetic equations for reactions 1 and 2 are as Eqs. (S20-S21)[10]:

$$r_1 = 6.3 \times 10^{10} \exp(-\frac{52000}{R\,T})\, c^D (c^A)^{0.5}\,, \tag{S20}$$

$$r_2 = 3.0 \times 10^9 \exp(-\frac{52000}{R\,T})(c^C)^2. \tag{S21}$$

Here, the reaction rates $r$ are expressed in kmol m$^{-3}$ s$^{-1}$ and the concentrations $c$ are in kmol m$^{-3}$. Temperature is in K. $A$, $C$ and $D$ represent $H_2$, $C_6H_6$ and $C_7H_8$ respectively. Furthermore, $B$ and $E$ will be used to represent $CH_4$ and $C_{10}H_{12}$ respectively.

### 6.2 Superstructure for the HDA process with a differential reactor model

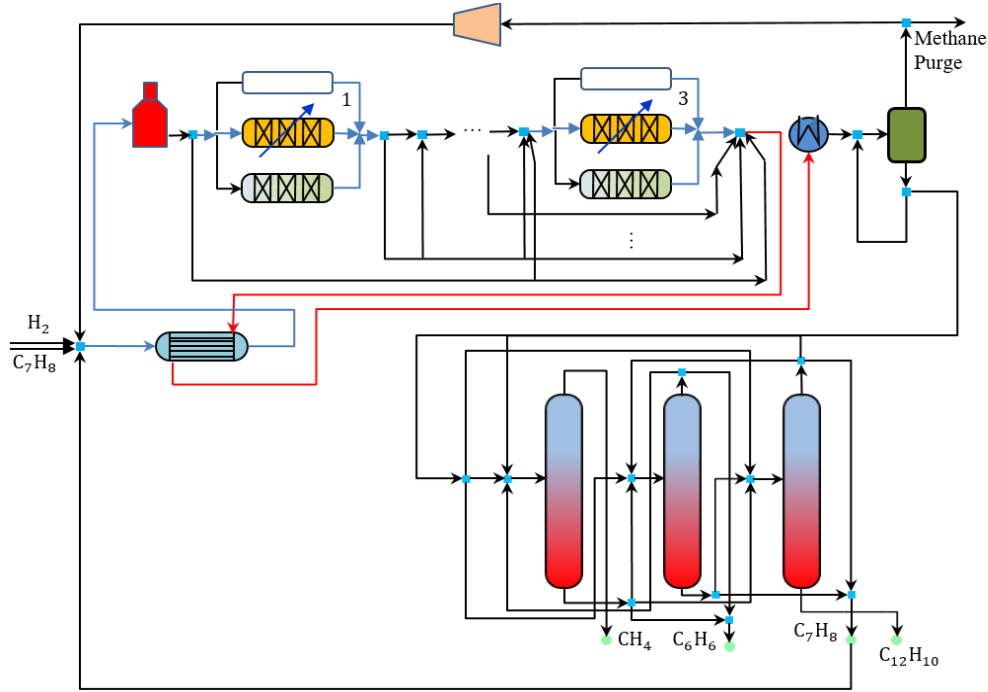

Figure S9 Superstructure for the synthesis of HDA process using differential reactor model

### 6.3 Orthogonal collocation finite element method model for PFR reactors

The collocation equation for mass balance is

$$\sum_{j=0}^{K} f_{m,j}^i \frac{d\ell_j(\tau_k)}{d\tau} = V\,d \sum_{p=1}^{2} \nu_{p,i} r_{m,k,p}$$

$$\forall m \in \mathbf{M}, k \in \mathbf{K} \backslash \{0\}, i \in \mathbf{I} \tag{S22}$$

where $\mathbf{M} := \{1, 2, \cdots M\}$ denotes the set of elements used to divide the reactor along the length direction, and the length of each element is $d = \frac{1}{M}$. $\mathbf{K} := \{0, 1, 2 \cdots K\}$ denotes the set of collocation points in each element, and $\mathbf{I} := \{A, B, C, D, E\}$ represents the set of components. $V$ (m$^3$) is the volume of the reactor. $f_{m,j}^i$ (kmol m$^{-3}$) denotes the molar flow rate of component $i$ at collocation point $j$ of the element $m$. $\nu_{p,i}$ is the stoichiometric coefficient shown in Table S16.

Table S16 Stoichiometric coefficients for HDA reactions

| $p$ | $i$ | $A$ | $B$ | $C$ | $D$ | $E$ |
|---|---|---|---|---|---|---|
| 1 | | -1 | 1 | 1 | -1 | 0 |
| 2 | | 1 | 0 | -2 | 0 | 1 |

$0 \leq \tau_k \leq 1$ is the location of the collocation point $k$. $\ell_j(\tau)$ represents the Lagrange polynomial

$$\ell_j(\tau) = \prod_{k=0, \neq j}^{K} \frac{\tau - \tau_k}{\tau_j - \tau_k}, \qquad\qquad j \in \mathbf{K} \qquad\qquad (S23)$$

and $\frac{d\ell_j(\tau_k)}{d\tau} = \frac{d\ell_j(\tau)}{d\tau}|_{\tau = \tau_k}$.

The continuity equations for the component flow rate profile are written as

$$f_{m+1,0}^i = \sum_{j=0}^{K} \ell_j(1) f_{m,j}^i, \qquad\qquad \forall m \in \mathbf{M} \backslash \{M\}, i \in \mathbf{I} \qquad (S24)$$

$$f^{out,i} = \sum_{j=0}^{K} \ell_j(1) f_{M,j}^i, \ \ f_{1,0}^i = f^{in,i}, \qquad \forall i \in \mathbf{I} \qquad\qquad (S25)$$

where $f^{out,i}$ and $f^{in,i}$ are the component flow rates in the reactor inlet and outlet streams respectively.

Similarly, the collocation equation and continuity equations for the enthalpy balance are

$$\sum_{j=0}^{K} H_{m,j} \frac{d\ell_j(\tau_k)}{d\tau} = V \, hq_{m,k}, \qquad\qquad \forall m \in \mathbf{M}, k \in \mathbf{K} \backslash \{0\} \qquad (S26)$$

$$H_{m+1,0} = \sum_{j=0}^{K} \ell_j(1) H_{m,j}, \qquad\qquad \forall m \in \mathbf{M} \backslash \{M\} \qquad (S27)$$

$$H^{out} = \sum_{j=0}^{K} \ell_j(1) H_{M,j}, \ \ H_{1,0} = H^{in}, \qquad\qquad\qquad (S28)$$

where $H_{m,j}$ denotes the enthalpy flow rate (kW) at collocation point $j$ of the element $m$, and $H^{in}, H^{out}$ represent the enthalpy flow rates of the reactor inlet and outlet streams. $q_{m,k}$ (kW m$^{-1}$) is the input energy from utility along the length direction. $H_{m,k}, H^{in}, H^{out}$ are calculated by

$$H_{m,k} = \sum_{i \in \mathbf{I}} f_{m,j}^i h_{m,j}^i, \qquad\qquad \forall m \in M, k \in K \qquad\qquad (S29)$$

$$H^{in} = \sum_{i \in \mathbf{I}} f^{in,i} h^{in,i}, \ H^{out} = \sum_{i \in \mathbf{I}} f^{out,i} h^{out,i}. \qquad\qquad (S30)$$

Here, $h_{m,j}^i$, $h^{in,i}$, $h^{out,i}$ are specific enthalpy of a component, which can be got from the Eq. (S17) with corresponding temperatures $T_{m,j}$, $T^{in}$ and $T^{out}$ applied.

The total energy input Q (kW) to the reactor can also be got through applying the collocation equation and continuity equations, as follows:

$$\sum_{j=0}^{K} Q_{m,j} \frac{d\ell_j(\tau_k)}{d\tau} = V\, h q_{m,k}, \qquad\qquad \forall m \in \mathbf{M}, k \in \mathbf{K}\backslash\{0\} \tag{S31}$$

$$Q_{m+1,0} = \sum_{j=0}^{K} \ell_j(1) Q_{m,j}, \qquad\qquad \forall m \in \mathbf{M}\backslash\{M\} \tag{S32}$$

$$Q^{out} = \sum_{j=0}^{K} \ell_j(1) Q_{M,j}, \;\; Q_{1,0} = 0\;, \tag{S33}$$

Here, $Q_{m,j}$ is the energy input from the reactor inlet to the collocation point $j$ of the element $m$. $Q^{out}$ is the total energy input to the reactor.

The reaction rate is got by

$$r_{m,k,1} = 6.3 \times 10^{10} \exp(-\frac{52000}{R\, T_{m,k}})\, c_{m,k}^D (c_{m,k}^A)^{0.5}\;, \;\; \forall m \in \mathbf{M}, k \in \mathbf{K} \tag{S34}$$

$$r_{m,k,2} = 3.0 \times 10^{9} \exp(-\frac{52000}{R\, T_{m,k}})(c_{m,k}^C)^2, \qquad \forall m \in \mathbf{M}, k \in \mathbf{K} \tag{S35}$$

where the component concentrations are from

$$c_{m,k}^i = \frac{f_{m,k}^i}{\sum_{i \in \mathbf{I}} f_{m,k}^i} \frac{P}{R\, T_{m,k}}. \qquad\qquad \forall m \in \mathbf{M}, k \in \mathbf{K}, i \in \mathbf{I} \tag{S36}$$

Here, $P$ is the operating pressure of the reactor.

Simplifications can be applied for given reactor types. For adiabatic reactors, $q_{m,k}$, $Q_{m,k}$ and $Q_{out}$ are all fixed at 0. For isothermal reactors, $T_{m,k} = T$.

## 6.4 Computational results

Table S17 Comparative results from BB, HCBb-FP, and HCBB-RB for Example 4

| Result | BB | HCBB-FP | HCBB-RB |
|---|---|---|---|
| | | Init 1 | |
| Profit | 4.904 | 4.904 | 4.904 |
| Profit after rounding | 4.903 | 4.903 | 4.903 |
| Relative change | 0.02% | 0.02% | -0.02% |
| | | Init 2 | |
| Profit | 4.904 | 4.904 | 4.904 |
| Profit after rounding | 4.903 | 4.903 | 4.903 |
| Relative change | 0.02% | 0.02% | 0.02% |
| | | Init 3 | |
| Profit | 4.904 | 4.904 | 4.904 |
| Profit after rounding | 4.903 | 4.903 | 4.903 |
| Relative change | 0.02% | 0.02% | 0.02% |
| | | Init 4 | |
| Profit | 4.904 | 4.904 | 4.904 |
| Profit after rounding | 4.903 | 4.903 | 4.903 |
| Relative change | 0.02% | 0.02% | 0.02% |
| | Init 5 | | |
| Profit | 4.901 | 4.901 | 4.901 |
| Profit after rounding | 4.899 | 4.899 | 4.899 |
| Relative change | 0.04% | 0.04% | 0.04% |
| | Init 6 | | |
| Profit | 4.309 | 4.903 | 4.904 |
| Profit after rounding | 4.306 | 4.902 | 4.903 |
| Relative change | 0.1% | 0.02% | 0.02% |

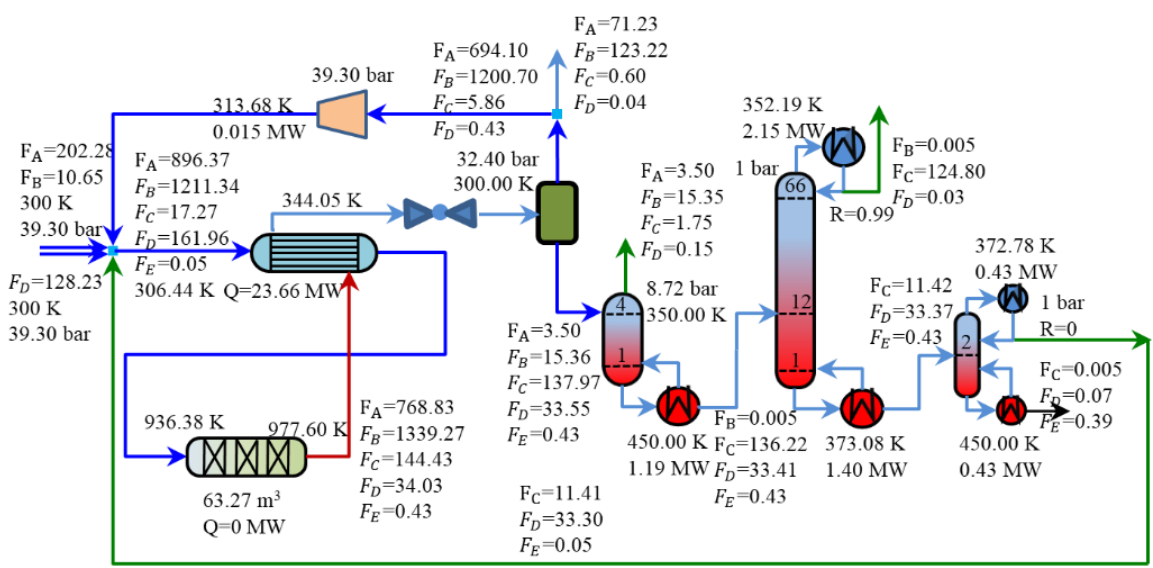

Figure S10 Optimal design of HDA process with differential reactor model where the unit of flow rate is kmol h$^{-1}$ and A, B, C, D and E represent $H_2$, $CH_4$, $C_6H_6$, $C_7H_8$ and $C_{12}H_{10}$

respectively.